\documentstyle[11pt,amstex,righttag,theorem,amssymb]{article}
\newif\ifafour
\afourtrue

\makeatletter
\ifafour\vsize=9.4truein
\if@twoside
\else
\fi
\columnwidth \textwidth
\columnseprule \z@
\skip\footins 10pt plus 2pt minus 4pt
\else
\if@twoside
\else
\fi
\fi
\def\@plus{plus}
\def\@string{\expandafter\@gobble\string}
\def\@oparg#1[#2]{\@ifnextchar[{#1}{#1[#2]}}
\long\def\@ifempty#1{%
 \expandafter\ifx\@car#1@\@nil @\empty
  \expandafter\@leftmark\else\expandafter\@rightmark\fi}
\long\def\@ifnotempty#1{\@ifempty{#1}{}}
\long\def\@leftmark#1#2{#1}
\long\def\@rightmark#1#2{#2}
\def\@andify#1#2{%
  \expandafter\in@\expandafter\and\expandafter{#2}%
  \ifin@
  \begingroup\global\toks\@ne{}
  \def\and##1\and{\@ifnotempty{##1}{\advance\@tempcnta\@ne\and}}%
  \@tempcnta\z@\expandafter\and#2\and%
  \ifnum\@tempcnta>\tw@
    \def\and##1\and{\advance\@tempcnta\m@ne
      \global\toks\@ne\expandafter{\the\toks\@ne ##1}%
      \edef\@tempa{\the\toks\@ne\ifnum\@tempcnta=\@ne\unskip,#1\else
          \ifnum\@tempcnta>\@ne\unskip,\ \fi\fi}%
      \global\toks\@ne\expandafter{\@tempa \ignorespaces}%
      \ifnum\@tempcnta>\z@\expandafter\and\fi}%
  \else 
    \def\and##1\and{\advance\@tempcnta\m@ne
       \global\toks\@ne\expandafter{\the\toks\@ne ##1}%
       \edef\@tempa{\the\toks\@ne
         \ifnum\@tempcnta=\@ne\unskip#1\ignorespaces\fi}%
       \global\toks\@ne\expandafter{\@tempa}%
       \ifnum\@tempcnta>\z@\expandafter\and\fi}%
  \fi 
 \expandafter\and#2\relax\endgroup
 \edef#2{\the\toks\@ne}%
\fi 
}
\def\uppercasetext@#1{%
   {\spaceskip1.3\fontdimen2\the\font plus1.3\fontdimen3\the\font
      \skipmath@#1$\skipmath@$}}

\def\upchars@{\def\ss{SS}\let\i=I\let\j=J\let\ae\AE\let\oe\OE
  \let\o\O\let\aa\AA\let\l\L}

\def\skipmath@#1$#2${\skipmath@b#1\(\skipmath@b\)%
  \ifx\skipmath@#2\else$#2$\expandafter\skipmath@\fi}

\def\skipmath@b#1\(#2\){\uppercase{#1}%
  \ifx\skipmath@b#2\else\(#2\)\expandafter\skipmath@b\fi}
\def\today{\ifcase\month\or
 January\or February\or March\or April\or May\or June\or
 July\or August\or September\or October\or November\or December\fi
 \space\number\day, \number\year}
\def\LaTeX{L\raise.42ex\hbox{\kern-.33em\protect\scriptsize A\kern-.15em}\TeX}
\def\@ifundefined#1{\expandafter\ifx\csname#1\endcsname\relax
        \expandafter\@leftmark
        \else\expandafter\@rightmark\fi}
\def\defaultfont{\family\default@family \series\default@series
  \shape\default@shape \selectfont}
\def\title{\@dblarg\@xtitle}
\def\@title{}
\def\@xtitle[#1]#2{\def\sh@rttitle{\def\\{\unskip\space\ignorespaces}%
   \ignorespaces#1\unskip}%
   \def\@title{\ignorespaces#2\unskip}}

\def\sh@rttitle{}
\def\author{\@dblarg{\@xauthor}}
\def\@authors{}
\def\@xauthor[#1]#2{\expandafter\def\expandafter
  \@authors\expandafter{\@authors#2\and}%
  \@ifnotempty{#1}{\expandafter\def\expandafter\sh@rtauthor
       \expandafter{\sh@rtauthor#1\and}}%
}
\def\sh@rtauthor{}
\newtoks\@addresstoks
\def\address{\@oparg\@xaddress[]}
\def\@xaddress[#1]#2{%
  \@addresstoks\expandafter{\the\@addresstoks\address{#1}{#2}{}{}{}{}}}
\def\@setaddress{%
  \par\expandafter\@ifnotempty\expandafter{\the\@addresstoks}{%
\begingroup \small
  \def\\{\unskip, \ignorespaces}%
  \interlinepenalty100
  \def\address##1##2##3##4##5##6{%
\vfil\penalty-300\vfilneg
\bigskip\indent
    \@ifnotempty{##1}{(\ignorespaces##1\unskip) }%
    {\sc\ignorespaces##2}\par
    \@ifnotempty{##4}{\nobreak\indent{\it Current address}%
      \@ifnotempty{##3}{, \ignorespaces##3\unskip}\/: 
      ##4\par}%
    \@ifnotempty{##6}{\nobreak\indent{\it E-mail address}%
      \@ifnotempty{##5}{, \ignorespaces##5\unskip}\/: 
      ##6\par}%
  }
  \the\@addresstoks
  \endgroup
}
}
\def\curraddr{\@oparg\@xcurraddr[]}
\toksdef\adthree3
\def\@xcurraddr[#1]#2{%
  \expandafter\@ifempty\expandafter{\the\@addresstoks}%
  {\errmessage{\noexpand\curraddr can't precede \string\address.}}%
  {\begingroup \global\toks3{}\global\toks1{}%
    \def\address##1##2##3##4##5##6{\global\toks3\expandafter{%
      \the\expandafter\adthree\the\toks1 }%
      \global\toks1{\address{##1}{##2}{##3}{##4}{##5}{##6}}}%
    \the\@addresstoks
    \def\address##1##2##3##4##5##6{\global\toks3\expandafter{%
      \the\toks3 \address{##1}{##2}{#1}{#2}{##5}{##6}}}%
    \the\toks\@ne
    \endgroup \@addresstoks\toks\thr@@
}}

\def\email{\@oparg\@xemail[]}

\def\@xemail[#1]#2{%
  \expandafter\@ifempty\expandafter{\the\@addresstoks}%
  {\errmessage{\noexpand\email can't precede \string\address.}}%
  {\begingroup \global\toks1{}\global\toks3{}%
    \def\address##1##2##3##4##5##6{\global\toks3\expandafter{%
      \the\expandafter\adthree\the\toks1 }%
      \global\toks1{\address{##1}{##2}{##3}{##4}{##5}{##6}}}%
    \the\@addresstoks
    \def\address##1##2##3##4##5##6{\global\toks3\expandafter{%
      \the\toks3 \address{##1}{##2}{##3}{##4}{#1}{#2}}}%
    \the\toks\@ne
    \endgroup \@addresstoks\toks\thr@@
}}

\def\date#1{\def\@date{#1}}
\def\@date{}

\def\thanks#1{\ifx\@empty\@thanks
   \def\@thanks{#1}%
 \else \expandafter\def\expandafter\@thanks\expandafter
   {\@thanks\@@par#1}\fi
}
\def\@thanks{}

\def\dedicatory#1{\def\@dedicatory{#1}}
\def\@dedicatory{}

\def\keywords#1{\def\@keywords{#1}}
\def\@keywords{}

\def\subjclass#1{\def\@subjclass{#1}}
\def\@subjclass{}

\def\translator#1{\def\@translator{#1}}
\def\@translator{}

\def\footnoterule{\kern-.4\p@
        \hrule\@width 5pc\kern11\p@\kern-\footnotesep}

\def\@makefnmark{\hbox{$\m@th^{\@thefnmark}$}}

\def\@makefntext{\indent\@makefnmark}
\def\sloppy{\tolerance9999\relax}
\long\def\@footnotetext#1{\insert\footins{%
    \defaultfont\footnotesize
    \interlinepenalty\interfootnotelinepenalty
    \splittopskip\footnotesep \splitmaxdepth \dp\strutbox
    \floatingpenalty\@MM \hsize\columnwidth \sloppy
  \edef\@currentlabel{\p@footnote\@thefnmark}\@makefntext
 {\rule\z@\footnotesep\ignorespaces#1\unskip\strut\par}}}

\def\maketitle{\par
  \@topnum\z@ 
  \ifx\@empty\sh@rtauthor \let\sh@rtauthor\sh@rttitle\fi
  \begingroup
  \@maketitle
  \endgroup
  \@andify{ AND }\sh@rtauthor
  \thispagestyle{plain}%
  \c@footnote\z@
  \def\do##1{\let##1\relax}%
  \do\maketitle \do\@maketitle
  \do\title \do\@xtitle \do\@title
  \do\author \do\@xauthor \do\@authors
  \do\address \do\@xaddress
  \do\email \do\@xemail \do\curraddr \do\@xcurraddr
  \do\dedicatory \do\@dedicatory
  \do\thanks \do\@thanks
  \do\keywords \do\@keywords
  \do\subjclass \do\@subjclass
  \do\@andify
}
\def\@maketitle{%
  \defaultfont\normalsize
  \let\@makefnmark\relax  \let\@thefnmark\relax
  \ifx\@empty\@subjclass\else
   \@footnotetext{1991 {\it Mathematics Subject
     Classification}.\enspace
        \@subjclass.}\fi
  \ifx\@empty\@keywords\else
   \@footnotetext{{\it Key words and phrases.}\enspace \@keywords.}\fi
\ifx\@empty\@thanks\else
   \@footnotetext{\@thanks}\fi
\topskip66\p@ 
  \vtop{\centering{\baselineskip20\p@\LARGE\@title\@@par}%
   \global\dimen@i\prevdepth}%
  \prevdepth\dimen@i
  \ifx\@empty\@authors
  \else
    \baselineskip32\p@
    \vtop{\@andify{ and }\@authors
      \centering{\Large\@authors\@@par}%
         \global\dimen@i\prevdepth}\relax
    \prevdepth\dimen@i
  \fi
  \ifx\@empty\@dedicatory
  \else
    \baselineskip18\p@
  \vtop{\centering{\small\it\@dedicatory\@@par}%
      \global\dimen@i\prevdepth}\prevdepth\dimen@i
  \fi
  \ifx\@empty\@date\else
  \baselineskip24\p@
    \vtop{\centering\@date\@@par
      \global\dimen@i\prevdepth}\prevdepth\dimen@i
  \fi
  \normalsize
  \vskip 2em
  } 

\def\abstractname{Abstract}
\def\abstract{%
  \global\let\abstract\relax
  \defaultfont\small  \skip@28\p@ \advance\skip@-\lastskip
  \advance\skip@-\baselineskip \vskip\skip@
  \moveright 3pc\vtop \bgroup
  \advance \hsize -6pc
  \trivlist \labelsep.5em\item[\hskip\labelsep
    {\sc\abstractname}.]\ignorespaces
}

\def\endabstract{\endtrivlist
  \global\let\endabstract\relax
  \global\dimen@i\prevdepth \egroup \prevdepth\dimen@i
  \skip@32\p@\@plus 14\p@ \advance\skip@-\baselineskip
  \vskip\skip@
}

\def\enddocument{\@checkend{document}\par
  \ifx\@empty\@translator \else
    \addvspace{6\p@\@plus9\p@}%
    \hbox to\columnwidth{\hss\defaultfont\normalsize
     Translated by \expandafter\uppercasetext@\expandafter{\@translator}}%
  \fi
\removelastskip
\@setaddress
 \clearpage\begingroup
 \if@filesw \immediate\closeout\@mainaux
 \def\global\@namedef##1##2{}\def\newlabel{\@testdef r}%
 \def\bibcite{\@testdef b}\@tempswafalse \makeatletter\input \jobname.aux
 \if@tempswa \@warning{Label(s) may have changed.  Rerun to get
 cross-references right}\fi\fi\endgroup\deadcycles\z@\@@end}
\makeatother

\title{Dressing Orbits of Harmonic Maps}
\author{F.E. Burstall}
\address[F.E. Burstall]{School of Mathematical Sciences\\
University of Bath \\
Bath BA2 7AY\\
England}
\thanks{First author partially supported by DOE grant DE-FG02
86ER25015, NSF grant DMS-9011083 and EEC contract SC1-0105-C.}
\email{f.e.burstall@@maths.bath.ac.uk}
\author{F. Pedit}
\address[F. Pedit]{Department of Mathematics \\
University of Massachusetts \\
Amherst, MA 01003\\USA}
\thanks{Second author supported by NSF grant DMS-9205293.}
\email{franz@@gang.umass.edu}

\date{}

\theoremstyle{plain}
\newtheorem{thm}{Theorem}[section]
\newtheorem{prop}[thm]{Proposition}
\newtheorem{lem}[thm]{Lemma}
\newtheorem{cor}[thm]{Corollary}

\theorembodyfont{\rm\normalshape}
\newtheorem{rem}{Remark} 
\newtheorem{defn}{Definition} 
\newtheorem{ex}{Example} 
\newtheorem{pf}{Proof}  
\newtheorem{notation}{Notation} 
 
\newenvironment{proof}{\begin{pf}}{\hspace*{\fill}$\square$\end{pf}}

\newcommand{\eqref}[1]{{\rm\normalshape(\ref{#1})}}

\newcounter{mynumi}
\newenvironment{mynum}{\begin{list}%
{{\rm\normalshape(}{\sl\roman{mynumi}\/}{\rm\normalshape)}}%
{\usecounter{mynumi}\setlength{\rightmargin}{\leftmargin}}}%
{\end{list}}

\newcommand{\R}{\Bbb{R}}
\newcommand{\C}{\Bbb{C}\,}
\newcommand{\N}{\Bbb{N}}
\renewcommand{\P}{\Bbb{P}}
\newcommand{\Z}{\Bbb{Z}}

\renewcommand{\b}{\frak{b}}
\newcommand{\centr}{\frak{z}_{A}}
\newcommand{\Comm}{Z}
\newcommand{\comm}{\cal{Z}}
\newcommand{\g}{\frak{g}}
\newcommand{\gc}{\g^{\C}}
\newcommand{\GC}{G^{\C}}

\newcommand{\k}{\frak{k}}
\newcommand{\kc}{\k^{\C}}
\newcommand{\KC}{K^{\C}}
\renewcommand{\O}{\cal{O}}
\renewcommand{\l}{\lambda}
\renewcommand{\L}{\Lambda}
\newcommand{\m}{\frak{m}}
\newcommand{\mc}{\m^{\C}}

\newcommand{\Harm}{\cal{H}}
\newcommand{\Ext}{\cal{E}}
\newcommand{\Gauge}{\cal{K}}
\newcommand{\SU}{\mathrm{SU}}
\newcommand{\SL}{\mathrm{SL}}
\newcommand{\SO}{\mathrm{SO}}
\newcommand{\Sl}{\frak{sl}}
\newcommand{\id}{\mathrm{I}}
\newcommand{\ad}{\operatorname{ad}}
\newcommand{\Ad}{\operatorname{Ad}}

\newcommand{\ext}{\mathrm{d}}
\newcommand{\half}{\tfrac{1}{2}}

\newcommand{\del}{\partial}
\newcommand{\Abar}{\overline{A}}
\newcommand{\Loop}{\Lambda^{\epsilon}G_{\tau}}
\newcommand{\Eloop}{\Lambda^{\epsilon}_{E}G_{\tau}}
\newcommand{\Iloop}{\Lambda^{\epsilon}_{I}G_{\tau}}
\newcommand{\IBloop}{\Lambda^{\epsilon}_{I,B}G_{\tau}}
\newcommand{\Hloop}{\Lambda_{\mathrm{hol}}G_{\tau}}
\newcommand{\gloop}{\Lambda^{\epsilon}\g_{\tau}}
\newcommand{\eloop}{\Lambda^{\epsilon}_{E}\g_{\tau}}
\newcommand{\iloop}{\Lambda^{\epsilon}_{I}\g_{\tau}}
\newcommand{\ibloop}{\Lambda^{\epsilon}_{I,\b}\g_{\tau}}

\newcommand{\initloop}{\Lambda^{\epsilon}_{-1,\infty}}
\newcommand{\rloop}{\Lambda\g_{\tau}}
\newcommand{\Orbit}{\cal{O}_{A}}
\newcommand{\dOrbit}{\cal{O}_{A}^{(d)}}
\newcommand{\Stab}{\Gamma^{\epsilon}_{I,B}}
\newcommand{\stab}{\cal{G}_{I,B}}
\newcommand{\zbar}{\bar{z}}
\numberwithin{equation}{section}

\begin{document}

\maketitle

\section*{Introduction}

At the heart of the modern theory of harmonic maps from a Riemann
surface to a Riemannian symmetric space is the observation that, in
this setting, the harmonic map equations have a zero curvature
representation \cite{Poh76,Uhl89,ZakSha78} and so correspond to loops
of flat connections.  This fact was first exploited in the
mathematical literature by Uhlenbeck in her study \cite{Uhl89} of
harmonic maps $\R^{2}\to G$ into a compact Lie group $G$.  Uhlenbeck
discovered that harmonic maps correspond to certain holomorphic maps,
the {\em extended solutions\/}, into the based loop group $\Omega G$
and used this to define an action of a certain loop group on the
space of harmonic maps.  However, the main focus of \cite{Uhl89} was
on harmonic maps of a $2$-sphere and for these maps the action
reduces to an action of a finite-dimensional quotient group (see also
\cite{BerGue91,GueOhn93}).

In another direction, the zero curvature representation has been
central to recent progress in the understanding of the harmonic map
equations as soliton equations, i.e. as completely integrable
Hamiltonian PDE . By solving certain Lax flows on loop algebras a
rather complete description of all harmonic tori in symmetric spaces
and Lie groups has been obtained. Of particular importance in this
approach are the harmonic maps of finite type: these arise from Lax
flows on finite dimensional subspaces of loop algebras and correspond
to linear flows on Jacobians of certain algebraic curves
\cite{BolPedWoo93,Bur93,BurFerPed93,FerPedPin92,Hit90,PinSte89}.
Among these harmonic maps, we further distinguish those of {\em
semisimple\/} finite type (see Section~\ref{sec:ft} below) which are
characterised by a semisimplicity condition on their derivative.
Semisimple finite type harmonic maps account for all non-conformal
harmonic tori in rank one symmetric spaces of compact type
\cite{BurFerPed93}, all non-isotropic harmonic tori in spheres and
complex projective spaces \cite{Bur93}, and all doubly periodic
solutions to the abelian affine Toda field equations for simple Lie
groups \cite{BolPedWoo93}.

The purpose of this paper is to describe some interactions between
these two approaches.  Our starting point is the fact that underlying
all of the above results is the existence of Iwasawa type
decompositions of the loop groups and algebras concerned.  On the one
hand, the Lax equations mentioned above arise from an Iwasawa
decomposition of certain twisted loop algebras via the
Adler--Kostant--Symes scheme \cite{BurPed94}.  On the other hand, the
loop group action of Uhlenbeck is essentially the {\em dressing\/}
action arising from the Iwasawa decompositions of the corresponding
loop groups \cite{GueOhn93}.  Moreover, a bridge beween these
constructions is provided by Symes's formula for the solution of the
Lax equations: in this construction, first applied by Symes
\cite{Sym80} to solve the open Toda lattice, projection of certain
complex geodesics on a factor in the Iwasawa decomposition yields
(extended framings of) harmonic maps.  This provides a map from a
certain subspace of the loop algebra to the space of harmonic maps
which intertwines the dressing and adjoint actions.

This set-up is very familiar in soliton theory where the dressing
method was first developed.  Here the idea is to use the dressing
action to construct new solutions from old and in many cases the
dressing orbits through trivial or {\em vacuum} solutions account for
all the solutions one is interested in. It is this theme that we
develop in the present paper.

We treat (primitive) harmonic maps $\R^{2}\to G/K$ where $G/K$ is a
($k$-)symmetric space.  In this context, our vacuum solutions are the
maps $f^{A}:\R^{2}\to G/K$ defined by
\[
f^{A}(z)=\exp(zA+\overline{zA})K,
\]
where $A$ is an element of the Lie algebra of $G$ satisfying
$[A,\overline{A}]=0$.  These maps (the geometry of which has recently
been studied by Jensen--Liao \cite{JenLia94} in the case $G/K=\C
P^{n}$) are equivariant with respect to actions of the abelian group
generated by $A$ and $\overline{A}$.

Our main results concern the orbit $\Orbit$ of $f^{A}$
under the dressing action.  This orbit is infinite-dimensional in
contrast to those studied by Uhlenbeck and we show that every
harmonic map of semisimple finite type lies in some $\Orbit$.
As a special case, we deduce that every non-isotropic harmonic torus
in a sphere or complex projective space is dressing equivalent to a
vacuum solution.

A distinctive feature of the orbits $\Orbit$ is that they admit a
hierarchy of commuting flows (conservation laws).  We show that this
hierarchy can be used to characterise the harmonic maps of finite
type: a harmonic map in $\Orbit$ is of finite type if and only if its
orbit under the hierarchy is finite-dimensional.

In all these results, just as in those of \cite{Bur93,BurFerPed93},
essential use is made of the semisimplicity assumption on the
derivative of the harmonic maps.  In an appendix, we examine the
situation when this assumption is dropped and discover intriguing
relationships between this case and harmonic maps of finite uniton
number in the sense of Uhlenbeck.

Special cases and partial versions of some of our results already
exist in the literature: non-conformal harmonic tori in
$S^{2}$---the Gauss maps of contant mean curvature tori---were
studied by Dorfmeister--Wu \cite{DorWu93} and a similar analysis of
non-superminimal minimal tori in $S^{4}$ was performed by Wu
\cite{Wu93}.  However, even in these cases, our methods are different
and we believe them to be more transparent.

Most of this research was carried out while both authors visited the
SFB 288 at TU-Berlin and the first author visited GANG at the
University of Massachusetts, Amherst.  We would like to thank the
members of both institutions for their support and hospitality.

Finally, the first author's understanding of the matters treated
herein was enhanced by conversations with Ian McIntosh and Martin
Guest for which he takes this opportunity to thank them.

\begin{notation}
Throughout this work, when a Lie group is denoted by an upper case
letter, its Lie algebra will be denoted by the corresponding lower
case gothic letter.  Thus $G$ is a Lie group with Lie algebra $\g$.
\end{notation}

\section{Primitive harmonic maps and their extended framings}
\label{sec:harm}

We are going to study primitive harmonic maps of a Riemann surface
into a $k$-symmetric space.  Such maps (or, rather, their framings)
have a zero-curvature representation and so give rise to maps into a
loop group.  This construction is fundamental for everything that
follows so we begin by reviewing this circle of ideas to establish
notation and to give a context for our results.

\subsection{Primitive harmonic maps}
\label{sec:prim}

Let $G$ be a compact semisimple Lie group.  A (regular) $k$-symmetric
$G$-space \cite{Kow80} is a coset space $N=G/K$ where
$(G^{\tau})_{0}\subset K\subset G^{\tau}$ for some automorphism
$\tau:G\to G$ of finite order $k\geq 2$.
\begin{ex}
A $2$-symmetric space is just a Riemannian symmetric space of compact
type.
\end{ex}
In general, the $k$-symmetric spaces form a large class of reductive
homogeneous spaces which include the generalised flag manifolds (that
is, $G/K$ where $K$ is the centraliser of a torus).

The automorphism $\tau$ induces a $\Z_{k}$-grading of $\gc$:
\begin{equation}
\gc=\sum_{\ell\in\Z_{k}}\g_{\ell},
\label{e:grad}
\end{equation}
where, setting $\omega=e^{2\pi i/k}$, $\g_{\ell}$ is the
$\omega^{\ell}$-eigenspace of (the derivative of) $\tau$.  We have
$\g_{0}=\kc$, $\overline{\g_{\ell}}=\g_{-\ell}$ and
\[
[\g_{j},\g_{\ell}]=\g_{j+\ell},
\]
where all arithmetic is modulo $k$.  In particular, defining
$\m\subset\g$ by $\mc=\sum_{\ell\in\Z_{k}\setminus\{0\}}\g_{\ell}$,
we have a reductive decomposition:
\begin{equation}
\g=\k\oplus\m.
\label{e:red}
\end{equation}
\begin{ex}
When $k=2$, $\g_{1}=\g_{-1}=\mc$ so that $[\m,\m]\subset\k$ and
\eqref{e:red} is the familiar symmetric decomposition.
\end{ex}

The decomposition \eqref{e:grad} induces a $G$-invariant
decomposition of the tangent bundle of $N$ which is non-trivial when
$k>2$.  Indeed, set $o=eK\in N$ and let $p=g\cdot o\in N$.  Then the
map $\g\to T_{p}N$ given by
\[
\xi\mapsto\left.\frac{d}{dt}\right|_{t=0}\exp t\xi\cdot p
\]
is a surjection with kernel $\Ad g\,\k$ and so gives an isomorphism
between $\Ad g\,\m\subset\g$ and $T_{p}N$.
\begin{notation}
If $\frak{l}\subset\g$ is an $\Ad K$-invariant subspace, we denote by
$[\frak{l}]$ the sub-bundle of $N\times\g$ with fibres given by
\[
[\frak{l}]_{g\cdot o}=\Ad g\,\frak{l}.
\]
\end{notation}
With this notation, we identify $[\m]$ with $TN$ and have a
decomposition
\[
TN^{\C}=\sum_{\ell\in\Z_{k}\setminus\{0\}}[\g_{\ell}].
\]

\begin{defn}
\cite{Bur93} A map $f:M\to N$ of a Riemann surface into a
$k$-symmetric space is {\em primitive\/} if $\ext
f(T^{1,0}M)\subset[\g_{-1}]$.
\end{defn}
Note that when $k=2$, $[\g_{-1}]=TN^{\C}$ and the primitivity
condition is vacuous.

The study of primitive maps can be motivated by the following
considerations: firstly, they arise naturally as twistor lifts
(prolongations) of certain harmonic maps into Riemannian symmetric
spaces \cite{Bur93,BurPed94}.  Secondly, there is a close
relationship between primitive maps and solutions of the affine Toda
field equations, both abelian and non-abelian.  In particular, there
is an essentially bijective correspondence between certain primitive
maps into the full flag manifold $G/T$ (modulo the left action of
$G$) and solutions to the abelian affine Toda field equations
\cite{BolPedWoo93,BurPed94} (in this case, $\tau$ is the
Coxeter--Killing automorphism).

In \cite{BurPed94}, it is shown that primitive maps are harmonic with
respect to suitable invariant metrics on $G/K$.  Moreover, it is
shown that the structure equations for such maps have the same form
as those for harmonic maps into $2$-symmetric spaces.  This motivates
the following definition:
\begin{defn}
\cite{BurPed94} A map $f:M\to N$ of a Riemann surface into a
$k$-symmetric space is {\em primitive harmonic\/} if $k=2$ and $f$ is
harmonic or $k>2$ and $f$ is primitive.
\end{defn}

\subsection{Extended framings}
\label{sec:fram}

Henceforth, we assume that the Riemann surface $M$ is contractible
(in our applications, we shall take $M=\R^{2}$).  As a consequence,
all maps $M\to G/K$ have global framings $g:M\to G$.  We will study
primitive harmonic maps via their framings.

Let $\pi:G\to G/K$ be the coset projection and $f:M\to G/K$ a
primitive harmonic map with framing $g:M\to G$, thus $f=\pi\circ g$.
Let $\alpha=g^{-1}\ext g$ be the pull-back by $g$ of the
Maurer--Cartan form of $G$.  We have a decomposition of $\alpha$
according to the eigenspace decomposition \eqref{e:grad} of $\gc$:
\[
\alpha=\sum_{\ell\in\Z_{k}}\alpha_{\ell}
\]
which, for a primitive harmonic map, reduces to
\[
\alpha=\alpha'_{-1}+\alpha_{0}+\alpha''_{1},
\]
where $\alpha'_{-1}$ is a $\g_{-1}$-valued $(1,0)$-form and
$\alpha''_{1}=\overline{\alpha'_{-1}}$ \cite{BurPed94}.  Moreover,
the condition that $f$ be primitive harmonic amounts to demanding
that the Maurer--Cartan equation for $\alpha$
\begin{equation}
\ext\alpha+\half[\alpha\wedge\alpha]=0
\label{e:mc}
\end{equation}
decouples into three equations:
\begin{gather}\label{e:phall}
\ext\alpha'_{-1}+[\alpha_{0}\wedge\alpha'_{-1}]=0
\tag{$\theequation'_{\m}$}\label{e:phm'}\\
\ext\alpha_{0}+\half[\alpha_{0}\wedge\alpha_{0}]+
[\alpha'_{-1}\wedge\alpha''_{1}]=0
\tag{$\theequation_{\k}$}\label{e:phk}\\
\ext\alpha''_{1}+[\alpha_{0}\wedge\alpha''_{1}]=0.
\tag{$\theequation''_{\m}$}\label{e:phm''}
\end{gather}
(For $k=2$, the equations \eqref{e:phm'} and \eqref{e:phm''} are the
harmonic map equations while, for $k>2$, they are the projections of
\eqref{e:mc} onto $\g_{-1}$ and $\g_{1}$.)

For $\lambda\in\C^{*}$, define a $\gc$-valued $1$-form by
\[
\alpha_{\lambda}=\lambda^{-1}\alpha'_{-1}+\alpha_{0}+\lambda\alpha''_{1}.
\]
By construction, $\alpha_{\lambda}$ enjoys the following properties:
\begin{enumerate}
\item $\alpha_{\lambda=1}=\alpha$.
\item For all $\lambda\in\C^{*}$,
$\tau\alpha_{\lambda}=\alpha_{\omega\lambda}$.
\item For all $\lambda\in\C^{*}$,
$\overline{\alpha_{\lambda}}=\alpha_{1/\bar{\lambda}}$ so that
$\alpha_{\lambda}$ is $\g$-valued when $\lambda\in S^{1}$.
\end{enumerate}
The crucial observation is that $\alpha$ satisfies the equations
\eqref{e:phall} if and only if
\[
\ext\alpha_{\lambda}+\half[\alpha_{\lambda}\wedge\alpha_{\lambda}]=0,
\]
for all $\lambda\in\C^{*}$.  Thus, for each $\lambda$, we can
integrate the Maurer--Cartan equations to obtain a map
$F_{\lambda}:M\to\GC$, unique up to left translation by a constant,
satisfying $F_{\lambda}^{-1}\ext F_{\lambda}=\alpha_{\lambda}$.
Moreover, in view of the properties of $\alpha_{\lambda}$ listed
above, we may choose the constants of integration to ensure
\begin{enumerate}
\item $F_{1}=g$.
\item For all $\lambda\in\C^{*}$, $\tau
F_{\lambda}=F_{\omega\lambda}$.
\item For all $\lambda\in\C^{*}$,
$\overline{F_{\lambda}}=F_{1/\bar{\lambda}}$ where conjugation is the
Cartan involution of $\GC$ fixing $G$.  In particular,
$F_{\lambda}:M\to G$ when $\lambda\in S^{1}$.
\item For each $p\in M$, $\lambda\mapsto F_{\lambda}(p)$ is
holomorphic on $\C^{*}$.
\end{enumerate}
Otherwise said, we have defined a map $F$ from $M$ into the group
$\Hloop$ given by

\[
\Hloop=\{g:\C^{*}\to\GC\colon\text{$g$ is holomorphic and
$\overline{g(\lambda)}=g(1/\bar{\lambda})$, $g(\omega\lambda)=\tau
g(\lambda)$}\}.
\]
This prompts the following definition:
\begin{defn}
A map $F:M\to\Hloop$ is an {\em extended framing\/} if
\[
F^{-1}\ext F=\lambda^{-1}\alpha'_{-1}+\alpha_{0}+\lambda\alpha''_{1},
\]
with $\alpha'_{-1}$ a $(1,0)$-form on $M$, or equivalently, if
$\lambda F^{-1}\del F$ is holomorphic at $\lambda=0$ (here
$F^{-1}\del F$ is the $(1,0)$-part of $F^{-1}\ext F$).
\end{defn}
We have seen that any primitive harmonic map $f$ admits an extended
framing $F$ such that $F_{1}$ is a framing of $f$.  Conversely, it is
clear that when $F$ is an extended framing then $F_{1}$ frames a
primitive harmonic map.  (In fact, $F_{\lambda}$ frames a primitive
harmonic map for each $\lambda\in S^{1}$).

This correspondence between primitive harmonic maps and extended
framings can be made bijective modulo gauge transformations by
imposing base-point conditions.  Fix $p_{0}\in M$ and let
\[
\Harm=\{f:M\to G/K\colon\text{$f$ is primitive harmonic with
$f(p_{0})=o$}\}
\]
be the space of based primitive harmonic maps.  Similarly, let
\[
\Ext=\{F:M\to\Hloop\colon\text{$F$ is an extended framing with
$F(p_{0})\in K$}\}
\]
be the space of based extended framings (here we identify $K$ with
the constant elements of $\Hloop$).  It is then straightforward to
see that we have a bijective correspondence
\[
\Harm\cong\Ext/\Gauge
\]
where the gauge group $\Gauge=\mathrm{C}^{\infty}(M,K)$ acts by
point-wise multiplication on the right.

With all this in place, our constructions involving primitive
harmonic maps will be made at the level of the corresponding extended
framings.  Before turning to this, however, we pause to briefly
describe a class of primitive harmonic maps which will be important
in the sequel.

\subsection{Primitive harmonic maps of finite type}
\label{sec:ft}

In \cite{BurPed94}, building on earlier work of several authors
\cite{BolPedWoo93,Bur93,BurFerPed93,FerPedPin92,PinSte89}, we
described a method for constructing primitive harmonic maps
$\R^{2}\to G/K$ from commuting Hamiltonian flows on loop algebras.

The loop algebra in question is $\rloop$ given by
\[
\rloop=\{\xi:S^{1}\to\g\colon\quad\tau\xi(\lambda)=\xi(\omega\lambda)\}.
\]
Any $\xi\in\rloop$ has a Fourier decomposition
\[
\xi=\sum_{n\in\Z}\lambda^{n}\xi_{n}
\]
and we distinguish the finite-dimensional subspaces
$\Lambda_{d}\subset\rloop$ given by
\[
\Lambda_{d}=\{\xi\in\rloop\colon\text{$\xi_{n}=0$ for $|n|>d$}\}.
\]

Now fix $d\equiv 1\mod k$.  A {\em polynomial Killing field\/} is a
map $\xi:\R^{2}\to\Lambda_{d}$ satisfying the Lax equation
\begin{equation}
\ext\xi=[\xi,(\lambda^{-1}\xi_{-d}+r(\xi_{1-d}))\,\ext z+
(\lambda\xi_{d}+\overline{r(\xi_{1-d})})\,\ext\zbar].
\label{e:lax}
\end{equation}
Here $z$ is the usual holomorphic co-ordinate on $\R^{2}$ and
$r:\g_{0}\to\g_{0}$ is a certain linear map constructed from an
Iwasawa decomposition of $\g_{0}$ (see \cite{BurPed94} for more
details).

Polynomial Killing fields simultaneously integrate a pair of
commuting Hamiltonian vector fields on $\Lambda_{d}$ and so there is
a unique such having any prescribed value at $z=0$.  A polynomial
Killing field gives rise to a primitive harmonic map because
$\alpha_{\lambda}$ defined by
\begin{equation}\label{e:xi2a}
\alpha_{\lambda}=(\lambda^{-1}\xi_{-d}+r(\xi_{1-d}))\,\ext z+
(\lambda\xi_{d}+\overline{r(\xi_{1-d})})\,\ext\zbar
\end{equation}
satisfies the Maurer--Cartan equations.  We may therefore integrate
to get an extended framing $F$ with $F_{\lambda}^{-1}\ext
F_{\lambda}=\alpha_{\lambda}$ and so a primitive harmonic map.  The
primitive harmonic maps that arise in this way are said to be of {\em
finite type\/}.

There is a necessary condition for a primitive harmonic map $f$ to be
of finite type: the Lax equation \eqref{e:lax} implies that
$\xi_{-d}:\R^{2}\to\g_{-1}$ must take values in a single
$\Ad\KC$-orbit so that, from \eqref{e:xi2a}, we see that, for some,
and hence every, framing of $f$, $\alpha'_{-1}(\del/\del z)$ must
take values in a single orbit also.  For doubly periodic primitive
harmonic maps (that is, those which cover a map of a torus), this
condition is almost sufficient:
\begin{thm}
{\normalshape\cite{BurPed94}} A doubly periodic primitive harmonic
map $\R^{2}\to G/K$ is of finite type if, for any framing,
$\alpha'_{-1}(\del/\del z)$ takes values in an $\Ad\KC$-orbit of {\em
semisimple\/} elements of $\g_{-1}$.
\end{thm}

In view of this, we make the following definition:
\begin{defn}
A primitive harmonic map $\R^{2}\to G/K$ is of semisimple finite type
if it admits a polynomial Killing field $\xi$ with $\xi_{-d}$
semisimple.
\end{defn}
The primitive harmonic maps of semisimple type include:
\begin{mynum}
\item All doubly periodic non-conformal harmonic maps of $\R^{2}$
into a rank one symmetric space \cite{BurFerPed93}.
\item Twistor lifts of any doubly periodic non-isotropic harmonic map
of $\R^{2}$ into a sphere or complex projective space \cite{Bur93}
(see also \cite{BolPedWoo93,FerPedPin92}).  Here $G/K$ is a flag
manifold in most cases.
\item Twistor lifts of certain doubly periodic harmonic maps of
$\R^{2}$ into low dimensional complex Grassmannians and quaternionic
projective spaces \cite{Uda94,Uda94a}.
\end{mynum}

\section{Dressing actions of loop groups}
We have seen that primitive harmonic maps correspond to extended
framings $M\to\Hloop$.  Various loop groups have non-trivial actions
on $\Hloop$ which induce actions on extended framings.  In this
section we define these groups and collect some elementary facts
concerning their generalised Birkhoff factorisations.  Using these
facts, we will be able to define these actions on extended framings
as well as provide a simple construction of extended framings via an
analogue of the Symes formula for solutions of the Toda lattice
\cite{Sym80}.

\subsection{Decompositions of loop groups}
\label{sec:decomp}

A $k$-symmetric space is defined by the following ingredients, which
we fix once and for all:
\begin{enumerate}
\item A compact semisimple group $G$.
\item An automorphism $\tau:G\to G$ of finite order $k\geq 2$ with
fixed set $K$.
\item The primitive $k$-th root of unity $\omega=e^{2\pi i/k}$.
\end{enumerate}
Moreover, we fix an Iwasawa decomposition of the reductive group
$\KC$:
\[
\KC=KB
\]
where $B$ is a solvable subgroup of $\KC$.  Thus any element
$k\in\KC$ can be uniquely written as a product
\[
k=k_{K}k_{B}
\]
with $k_{K}\in K$, $k_{B}\in B$.

The loop groups of interest to us are spaces of smooth maps from a
pair of circles in $\C$ to $\GC$ which are equivariant with respect
to $\tau$ and satisfy a reality condition.  To define them, we fix
$0<\epsilon<1$ and partition the Riemann sphere
$\P^{1}=\C\cup\{\infty\}$ as follows: let $C_{\epsilon}$ and
$C_{1/\epsilon}$ denote the circles of radius $\epsilon$ and
$1/\epsilon$ about $0\in\C$ and define open sets by
\[
I_{\epsilon}=\{\lambda\in\P^{1}\colon|\lambda|<\epsilon\},\quad
I_{1/\epsilon}=\{\lambda\in\P^{1}\colon|\lambda|>1/\epsilon\},\quad
E^{(\epsilon)}=\{\lambda\in\P^{1}\colon\epsilon<|\lambda|<1/\epsilon\}.
\]
Now put $I^{(\epsilon)}=I_{\epsilon}\cup I_{1/\epsilon}$ and
$C^{(\epsilon)}=C_{\epsilon}\cup C_{1/\epsilon}$ so that
$\P^{1}=I^{(\epsilon)}\cup C^{(\epsilon)}\cup E^{(\epsilon)}$.

We define the group of smooth maps $\Loop$ by
\[
\Loop=\{g: C^{(\epsilon)}\to G^{\C}\colon\text{$g(\omega
\lambda)=\tau g(\lambda)$, $\overline{g(\lambda)}=g
(1/\bar{\lambda})$, for all $\lambda\in C^{(\epsilon)}$}\}.
\]
Here, again, the conjugation is the Cartan involution of $\GC$ which
fixes $G$.

\begin{rem}
Observe that the reality condition
$\overline{g(\lambda)}=g(1/\bar{\lambda})$ implies that $g\in\Loop$
is completely determined by its values on $C_{\epsilon}$ so that we
have an isomorphism between $\Loop$ and the group of
$\tau$-equivariant maps $C_{\epsilon}\to\GC$.  In particular, $\Loop$
becomes in this way a {\em complex\/} Lie group.  In what follows, we
shall use this isomorphism to identify elements of $\Loop$ with their
restrictions to $C_{\epsilon}$.
\end{rem}

We now define some subgroups of $G$:
\begin{align*}
\Eloop&=\{g\in\Loop\colon\text{$g$ extends holomorphically to
$g:E^{(\epsilon)}\to\GC$}\},\\ \Iloop&=\{g\in\Loop\colon\text{$g$
extends holomorphically to $g:I^{(\epsilon)}\to\GC$}\}.
\end{align*}
By unique continuation, any element $g$ of these subgroups satisfies
the reality and equivariance conditions
\[
\overline{g(\lambda)}=g(1/\bar{\lambda}),\qquad g(\omega\lambda)=\tau
g(\lambda),
\]
for all $\lambda$ in its domain of definition.  In particular,
$g\in\Iloop$ has $g(0)\in\KC$ and we distinguish the subgroup
\[
\IBloop=\{g\in\Iloop\colon g(0)\in B\}.
\]

The main tool in all our constructions is the following Iwasawa type
decomposition for the complex loop groups $\Loop$ which is due to
McIntosh \cite{McI94}.
\begin{thm} \label{thm:McInt}
Multiplication $\Eloop\times\IBloop\to\Loop$ is a diffeomorphism
onto.
\end{thm}

\begin{rem}
The important fact here is the surjectivity of the multiplication
which amounts to the assertion that a certain Riemann--Hilbert
problem is always solvable.  That this is indeed the case is a
consequence of the reality conditions we have imposed on our loops.
\end{rem}

\begin{rem}
The limiting case of Theorem~\ref{thm:McInt} as $\epsilon\to 1$ is
the more familiar assertion that a loop $S^{1}\to\GC$ can be
factorised as a product of a loop $S^{1}\to G$ and a loop which has a
holomorphic extension to $\{|\lambda|<1\}$.  This result is due to
Pressley--Segal \cite{PreSeg86} and was extended to the twisted
setting by Dorfmeister--Pedit--Wu \cite{DorPedWu94}.
\end{rem}

As a consequence of Theorem~\ref{thm:McInt}, we have a diffeomorphism
\[
\Eloop\cong\Loop/\IBloop
\]
so that $\Loop$ and, particularly, $\Iloop$ acts on $\Eloop$.  To
describe this action explicitly, note that any $g\in\Loop$ has a
unique factorisation
\[
g=g_{E}g_{I}
\]
with $g_{E}\in\Eloop$ and $g_{I}\in\IBloop$.  The action of $\Iloop$
on $\Eloop$ is now given by
\begin{equation} \label{e:dress}
g\#_\epsilon h = (gh)_E,
\end{equation}
for $g\in\Iloop$, $h\in\Eloop$.  We call this action the {\em
dressing action\/} on $\Eloop$.

\begin{rem}
These ideas fit into the general framework of dressing actions on
Poisson--Lie groups described by Lu--Weinstein \cite{LuWei90}: there
is a Poisson--Lie structure on $\Eloop$ for which $\Loop$ is the
double group and $\IBloop$ is the dual group.  In this context, our
dressing action is precisely the right dressing action of the dual
group in the sense of Lu--Weinstein made into a left action in the
usual way.  We shall return to these matters elsewhere.
\end{rem}

In our applications to harmonic maps, we need to factor out the
constant loops $K$ in $\Eloop$.  This is compatible with the dressing
action:
\begin{lem} \label{lem:desc}
The dressing action \eqref{e:dress} descends to an action on
$\Eloop/K $.
\end{lem}

\begin{proof}
First observe that $K$ is stable under the dressing action: if
$g\in\Iloop$ with $g_{0}=g(0)$ and $k\in K$, we use the Iwasawa
decomposition of $\KC$ to write
\[
gk=(g_{0}k)_{K}\left((g_{0}k)_{B}(k^{-1}g_{0}^{-1}gk)\right)
\]
so that
\[
g\#_{\epsilon}k=(g_{0}k)_{K}\in K.
\]
Now, for $h\in\Eloop$, we have
\[
g\#_{\epsilon}(hk)=(ghk)_{E}=\bigl((gh)_{E}(gh)_{I}k\bigr)_{E}=
(g\#_{\epsilon}h)\bigl((gh)_{I}\#_{\epsilon}k\bigr)
\]
and $(gh)_{I}\#_{\epsilon}k=\bigl((gh)_{I}(0)k\bigr)_{K}\in K$.
\end{proof}

To see how these constructions vary with $\epsilon$, note that, for
$0<\epsilon<\epsilon'<1$, restriction of the holomorphic extensions
provides injections
\[
\Lambda^{\epsilon'}_{I}G^{\C}_{\tau}\subset\Iloop,\qquad
\Eloop\subset\Lambda^{\epsilon'}_{E}G^{\C}_{\tau}
\]
and similarly, for $0<\epsilon<1$, we have
\[
\Hloop\subset\Eloop.
\]
Indeed, it is easy to see that
\begin{equation}
\Hloop=\bigcap_{0<\epsilon<1}\Eloop.
\label{e:holo-in-E}
\end{equation}
The dressing actions are compatible with these inclusions because the
following generalisation of a result of Guest--Ohnita
\cite{GueOhn93}:
\begin{prop} \label{prop:compat}
For $0<\epsilon<\epsilon'<1$,
$g\in\Lambda^{\epsilon'}_{I}G^{\C}_{\tau}\subset\Iloop$ and
$h\in\Eloop\subset\Lambda^{\epsilon'}_{E}G^{\C}_{\tau}$, we have
\[
g\#_{\epsilon'}h=g\#_{\epsilon}h\in\Eloop.
\]
\end{prop}
\begin{proof}
On $C^{(\epsilon')}$, we have
\[
g\#_{\epsilon'}h=gh(gh)^{-1}_{I^{(\epsilon')}}
\]
where the left hand side has a holomorphic extension to
$E^{(\epsilon')}$ while, since $h\in\Eloop$, the right hand side has
a holomorphic extension to $I^{(\epsilon')}\cap E^{(\epsilon)}$.  It
now follows from a theorem of Painlev\'e that $g\#_{\epsilon'}h$ has
a holomorphic extension to $E^{(\epsilon')}\cup
C^{(\epsilon')}\cup(I^{(\epsilon')}\cap
E^{(\epsilon)})=E^{(\epsilon)}$.  Thus $g\#_{\epsilon'}h\in\Eloop$
while $(gh)_{I^{(\epsilon')}}\in\IBloop$ and the proposition follows
from the uniqueness of the factorisation of $\Loop$.
\end{proof}
In particular, the action of $\Lambda_{I}^{\epsilon'}G^{\C}_{\tau}$
preserves each $\Eloop$ for $0<\epsilon<\epsilon'<1$ and so, taking
\eqref{e:holo-in-E} into account, we conclude
\begin{cor}\label{cor:holo-stable}
The action of each $\Eloop$ preserves $\Hloop$ and, for
$0<\epsilon<\epsilon'<1$,
$g\in\Lambda^{\epsilon'}_{I}G^{\C}_{\tau}\in\Iloop$ and $h\in\Hloop$,
we have
\[
g\#_{\epsilon'}h=g\#_{\epsilon}h.
\]
\end{cor}
\begin{notation}
Henceforth, we simply write $g\#h$ for the action on $\Hloop$.
\end{notation}
\begin{rem}
It follows from Corollary~\ref{cor:holo-stable} that we can take a
direct limit as $\epsilon\to0$ and so obtain an action on $\Hloop$ of
the group of germs at zero of $\tau$-equivariant maps $\C\to\GC$.
This action is very similar to the one discussed by Uhlenbeck
\cite{Uhl89} although she considers only the subgroup of rational
maps $\P^{1}\to\GC$ which satisfy the reality conditions and are
holomorphic at zero.
\end{rem}

Finally, we note that, by virtue of Lemma~\ref{lem:desc}, we have
\begin{cor}\label{cor:hol-desc}
The dressing action of each $\Iloop$ descends to an action on
$\Hloop/K$.
\end{cor}

\subsection{Symes formula: point-wise version}
\label{sec:symes1}

In \cite{Sym80}, Symes gave a formula for solutions of the Toda
lattice in terms of the projection of a one-parameter subgroup onto
one of the factors in an Iwasawa decomposition.  We shall see that a
similar formula holds for extended framings.  We begin by describing
the space of generators for these one-parameter subgroups.

The Lie algebras of the loop groups of the previous section are
simply the corresponding algebras of $\tau$-equivariant maps
$C^{(\epsilon)}\to\gc$ and are denoted $\gloop$, $\eloop$ and so on.
Throughout this section, we use the reality condition to identify
elements of $\gloop$ and its subalgebras with maps
$C_{\epsilon}\to\gc$.

Define the subspace $\initloop\subset\gloop$ by
\[
\initloop=\{\xi\in\gloop\colon\text{$\lambda\xi$ has a holomorphic
extension to $I_{\epsilon}$}\}.
\]
Thus $\initloop$ consists of those elements of $\gloop$ which extend
meromorphically to $I^{(\epsilon)}$ with at most simple poles at $0$
and $\infty$.

Observe that $\initloop$ is stable under the adjoint action of
$\Iloop$ on $\gloop$: if $g\in\Iloop$ and $\xi\in\initloop$ then
$\lambda\xi$ extends holomorphically to $I_{\epsilon}$ whence
\[
\lambda\Ad g(\xi)=\Ad g(\lambda\xi)
\]
does also.

Now define $\Phi_{\epsilon}:\initloop\to\Eloop$ by
\[
\Phi_{\epsilon}(\xi)=(\exp\xi)_{E},
\]
where the exponential map is defined point-wise:
$(\exp\xi)(\lambda)=\exp_{\GC}\xi(\lambda)$.

Arguing as in Proposition~\ref{prop:compat}, one proves
\begin{prop}
For $0<\epsilon<\epsilon'<1$ and
$\xi\in\Lambda^{\epsilon'}_{-1,\infty}\subset\initloop$, we have
\[
\Phi_{\epsilon'}(\xi)=\Phi_{\epsilon}(\xi)\in\Eloop.
\]
In particular, $\Phi_{\epsilon'}$ has image in
$\bigcap_{\epsilon<\epsilon'}\Eloop=\Hloop$.
\end{prop}
\begin{notation}
In view of this, we shall henceforth simply write $\Phi$ for
$\Phi_{\epsilon}$.
\end{notation}

A useful property of $\Phi$ is that it essentially intertwines the
adjoint and dressing actions of $\Iloop$:
\begin{prop}\label{prop:intw}
Let $g\in\Iloop$ and $\xi\in\initloop$.  Then
\[
\Phi(\Ad g\,\xi)=g\#\bigl(\Phi(\xi)k\bigr),
\]
where $k=\bigl((\exp\xi)_I(0)g(0)^{-1}\bigr)_K\in K$.
\end{prop}
Note that when $g(0)\in B$, $k=1$ so that $\Phi$ truly intertwines
the actions of $\IBloop$.
\begin{proof}
Using $\exp\Ad g\,\xi=g\exp(\xi)g^{-1}$, we have
\begin{align*}
\Phi(\Ad g\,\xi)&=(\exp\Ad g\,\xi)_{E}=(g\exp(\xi)g^{-1})_{E}\\
&=\bigl(g\Phi(\xi)(\exp\xi)_{I}g^{-1}\bigr)_{E}
=(g\Phi(\xi)k\widetilde{g})_{E}\\ &=g\#(\Phi(\xi)k)
\end{align*}
where $ \widetilde{g}\in \IBloop$ and $ k=((\exp
\xi)_I(0)g(0)^{-1})_K $ so that $ \exp(\xi)_I g^{-1}=k\widetilde{g}
$.
\end{proof}

In particular, denoting by $\overline{\Phi}:\initloop\to\Hloop/K$ the
composition of $ \Phi $ with the coset projection, we conclude from
Proposition~\ref{prop:intw} together with
Corollary~\ref{cor:hol-desc}:
\begin{cor}\label{cor:desc-intw}
$\overline{\Phi}$ intertwines the adjoint and (descended) dressing
actions of $\Iloop$:
\[
\overline{\Phi}(\Ad g\, \xi)=g\#\overline{\Phi}(\xi),
\]
for $\xi\in\initloop$, $g\in\Iloop$.
\end{cor}

\subsection{Dressing action on extended framings}
\label{sec:dress-ext}

The relevance of the results of Section~\ref{sec:decomp} to primitive
harmonic maps is that the point-wise dressing action of $\Iloop$ on
$\Hloop$ induces an action on extended framings that preserves gauge
orbits and the base-point condition.  We therefore arrive at an
action of $\Iloop$ on $\Harm$:
\begin{prop}\label{prop:dress-ext}

Let $g\in\Iloop$ and $F:M\to\Hloop$ be an extended framing.  Define
$g\#F:M\to\Hloop$ by
\[
(g\#F)(p)=g\#(F(p)),
\]
for $p\in M$.  Then
\begin{mynum}
\item $g\#F$ is also an extended solution.
\item If $F$ is based (that is, $F\in\Ext$) then so is $g\#F$.
\item If $k\in\Gauge$ then
\begin{equation}
g\#(Fk)=(g\#F)\widetilde{k}
\label{e:gauge}
\end{equation}
with $\widetilde{k}\in\Gauge$.
\end{mynum}
Thus $\Iloop$ acts on $\Harm=\Ext/\Gauge$.
\end{prop}
\begin{proof}
To see that $g\#F$ is an extended frame, write
\[
gF=ab
\]
where $a=g\#F$ and $b:M\to\IBloop$.  Then
\begin{align*}
a^{-1}\ext a&=\Ad b(F^{-1}\ext F-b^{-1}\ext b)\\ \intertext{so that}
\lambda a^{-1}\del a&=\Ad b(\lambda F^{-1}\del F-\lambda b^{-1}\del
b).
\end{align*}
Now all ingredients on the right hand side are holomorphic in
$\lambda$ on $I_{\epsilon}$ so that $\lambda a^{-1}\del a$ is also
whence $a$ is an extended framing.

Now suppose that $F\in\Ext$ so that $F(p_{0})\in K$.  In the proof of
Lemma~\ref{lem:desc} we saw that the action of $\Iloop$ preserves $K$
so that $(g\#F)(p_{0})=g\#(F(p_{0}))\in K$ as required.

Finally, \eqref{e:gauge} is an immediate consequence of
Lemma~\ref{lem:desc}.
\end{proof}
\begin{rem}
Proposition~\ref{prop:dress-ext} is a simple generalisation to our
setting of results of Uhlenbeck \cite{Uhl89} and Bergveldt--Guest
\cite{BerGue91} for harmonic maps into Lie groups.
\end{rem}

\subsection{Symes formula for extended framings}
\label{sec:symes-ext}

We now use the map $\Phi$ of Section~\ref{sec:symes1} to construct
based extended solutions $\R^{2}\to\Hloop$ from complex lines in
$\initloop$.  (In all that follows, we base maps of $\R^{2}$ at $0$.)

For $\eta\in\initloop$, define $F^{\eta}:\R^{2}\to\Hloop$ by
\[
F^{\eta}(z)=\Phi(z\eta)=(\exp z\eta)_{E}.
\]
\begin{prop}\label{prop:sym-xt}
$F^{\eta}$ is a based extended framing.
\end{prop}
\begin{proof}
We have
\[
\exp(z\eta)=F^{\eta}b,
\]
with $b:\R^{2}\to\IBloop$ so that
\begin{equation}
\label{e:sym-deriv}
\lambda(F^{\eta})^{-1}\del F^{\eta}= \Ad b(\lambda\eta\,\ext
z-\lambda b^{-1}\del b).
\end{equation}
Again all ingredients on the right hand side are holomorphic near
$\lambda=0$ (since $\eta\in\initloop$) so that
$\lambda(F^{\eta})^{-1}\del F^{\eta}$ is also whence $F^{\eta}$ is an
extended framing.  Finally, $F^{\eta}(0)=1\in K$ so that
$F^{\eta}\in\Ext$.
\end{proof}

Observe that primitive harmonic maps that arise from
Proposition~\ref{prop:sym-xt} in this way have the same restriction
on their framings as the maps of finite type: indeed, let
$\eta_{-1}=(\lambda\eta)(0)$ and put $b_{0}=b|_{\lambda=0}:\R^{2}\to
B$.  Now set $\lambda=0$ in \eqref{e:sym-deriv} to conclude that
\[
(F^{\eta})^{-1}\ext F^{\eta}=
\lambda^{-1}\alpha'_{-1}+\alpha_{0}+\lambda\alpha''_{1}
\]
with
\begin{equation}
\label{e:orbit}
\alpha'_{-1}=\Ad b_{0}(\eta_{-1})\,\ext z.
\end{equation}
Otherwise said, $\alpha'_{-1}(\del/\del z)$ takes values in a single
$B$-orbit.

In fact, all primitive harmonic maps of finite type arise from this
construction: let $f:\R^{2}\to G/K$ be a based primitive harmonic map
of finite type with polynomial Killing field
$\xi:\R^{2}\to\Lambda_{d}$ and let
$\eta=\lambda^{d-1}\xi(0)\in\initloop$.  The arguments of
\cite[Section~4.2]{BurPed94} carry over directly to our setting to
prove:
\begin{enumerate}
\item $F^{\eta}$ is an extended framing for $f$.
\item $\xi=\Ad(F^{\eta})^{-1}\xi(0)$.
\end{enumerate}
To summarise:
\begin{prop}\label{prop:ft}
$f:\R^{2}\to G/K$ is of finite type if and only if it admits an
extended framing of the form $F^{\eta}$ where, for some $d\equiv
1\mod k$, $\eta\in\lambda^{d-1}\Lambda_{d}\subset\initloop$
\end{prop}

Finally, we note the following useful fact: viewing the Symes formula as
a map
\[
\initloop\to\Harm,
\]
we see that this map intertwines the adjoint and dressing actions of
$\Iloop$.  Indeed, denoting by $[F^{\eta}]$ the gauge equivalence
class of $F^{\eta}$, Corollary~\ref{cor:desc-intw} gives
\begin{lem}\label{lem:intw-xt}
For $g\in\Iloop$ and $\eta\in\initloop$,
\[
[F^{\Ad g\,\eta}]=g\#[F^{\eta}].
\]
\end{lem}

\section{Dressing orbits of vacuum solutions}\label{sec:vac}

The main philosophy of the dressing construction
\cite{Wil85,ZakSha78} is to construct non-trivial solutions of
partial differential equations by applying a dressing action of a
suitable loop group to a trivial or {\em vacuum\/} solution.  In this
section, we define vacuum solutions for our problem and show that any
primitive harmonic map of semisimple finite type lies in the dressing
orbit of such a solution.

Some care must be exercised when developing the notion of a vacuum
solution of the primitive harmonic map equations: the simplest
harmonic maps are the constant maps but we learn from
Lemma~\ref{lem:desc} that the constant element of $\Harm$ is
preserved by the action of $\Iloop$ so that dressing will give us
nothing new.  (Similar remarks apply to the harmonic maps of finite
uniton number in the sense of Uhlenbeck \cite{Uhl89}, see the
appendix.)  However, if one is guided by the relationship between
primitive harmonic maps and Toda fields, one is led to the following
class of solutions which correspond to constant solutions of the
appropriate Toda field equations.
\begin{defn}
A {\em vacuum solution\/} is an extended framing of the form
$F^{\eta_{A}}$ where $\eta_{A}=\lambda^{-1}A$ with $A\in\g_{-1}$ and
$[A,\Abar]=0$.
\end{defn}
Notice that it follows from the vanishing of $[A,\Abar]$ that $A$ is
semisimple.

In this simple case, we can perform the Iwasawa decomposition
explicitly:
\[
\exp(z\eta_{A})=\exp(z\lambda^{-1}A+\zbar\lambda\Abar)\exp(-\zbar\lambda\Abar),
\]
on $C_{\epsilon}$ so that
\[
F^{\eta_{A}}(z)=\exp(z\lambda^{-1}A+\zbar\lambda\Abar)
\]
and the corresponding primitive harmonic map $f:\R^{2}\to G/K$
\[
f(z)=\exp(z\lambda^{-1}A+\zbar\lambda\Abar)o
\]
is equivariant for actions of the abelian group generated by $A$ and $\Abar$.

Denote by $\Orbit\subset\Harm$ the $\IBloop$-orbit of
$[F^{\eta_{A}}]$ so that
\[
\Orbit\cong\IBloop/\Stab
\]
where $\Stab$ is the stabiliser of $[F^{\eta_{A}}]$:
\[
\Stab=\{g\in\IBloop\colon\quad g\#[F^{\eta_{A}}]=[F^{\eta_{A}}]\}.
\]
\begin{rem}
Let us pause to justify our decision to study $\Orbit$ rather than
the slightly larger $\Iloop$-orbit of $[F^{\eta_{A}}]$.  We shall see
in Section~\ref{sec:flow} that $\Orbit$ admits a hierarchy of
commuting flows which can be used to characterise the primitive
harmonic maps of finite type.  These flows do not extend to the
$\Iloop$-orbit.  Moreover, since $\Iloop=K\IBloop$ and the dressing
action on $\Harm$ of $K$ is merely that induced by the action of $K$
on $G/K$, we see that we lose nothing essential by adopting this
approach.
\end{rem}
We will study $\Orbit$ by using the equivariance of the map
$\eta\to[F^{\eta}]:\initloop\to\Harm$ to replace the dressing action
of $\IBloop$ on $\Harm$ by the easier adjoint action on $\initloop$.
To accomplish this, we must describe the fibres of this map.  We
begin with the following observation:
\begin{lem}\label{lem:fibre}
Let $\zeta$, $\eta\in\initloop$.  Then $[F^{\zeta}]=[F^{\eta}]$ if
and only if $(\lambda\zeta)(0)=(\lambda\eta)(0)$ and
\begin{equation}
\label{e:comm}
(\ad\eta)^{n}\zeta\in\iloop,
\end{equation}
for all $n\geq1$.
\end{lem}
\begin{proof}
$[F^{\zeta}]=[F^{\eta}]$ if and only if $F^{\zeta}=F^{\eta}k$, for
some $k\in\Gauge$, and, using the definitions of $F^{\zeta}$ and
$F^{\eta}$, it is straight-forward to see that this is the case
precisely when
\[
e(z):=\exp(-z\zeta)\exp(z\eta)\in\Iloop,
\]
for $z\in\R^{2}$.  This, in turn, is the same as demanding that
\[
e^{-1}\ext e =(-\Ad\exp(-z\eta)\,\zeta+\eta)\,\ext z
\]
be $\iloop$-valued, that is,
\[
e^{-\ad z\eta}\zeta-\eta\in\iloop,
\]
for all $z\in\R^{2}$.  Expanding this last in powers of $z$ and
comparing coefficients proves the lemma.
\end{proof}

Applying this to the case where $\eta=\eta_{A}$ gives:
\begin{prop}\label{prop:fibreA}
$[F^{\zeta}]=[F^{\eta_{A}}]$ if and only if $(\lambda\zeta)(0)=A$ and
$[\zeta,A]=0$.
\end{prop}
\begin{proof}
Write $\zeta=\sum_{n\geq -1}\lambda^{n}\zeta_{n}$ on $C_{\epsilon}$.
Comparing coefficients of $\lambda$ in \eqref{e:comm} gives
\[
(\ad A)^{n}\zeta_{n-1}=0,
\]
for all $n\geq1$.  However, since $A$ is semisimple, $\ker(\ad
A)^{n}=\ker\ad A$ whence $[\zeta,A]=0$ as required.
\end{proof}

Taking Lemma~\ref{lem:intw-xt} into account, we obtain the following
characterisation of the fibres of $\initloop\to\Harm$ over $\Orbit$:
\begin{prop}\label{prop:fibreB}
For $\zeta\in\initloop$, $g\in\IBloop$,
$[F^{\zeta}]=g\#[F^{\eta_{A}}]\in\Orbit$ if and only if
$(\lambda\zeta)(0)=\Ad g(0)\,A$ and
\[
[\zeta,\Ad g\,A]=0.
\]
\end{prop}
\begin{proof}
By Lemma~\ref{lem:intw-xt}, $[F^{\zeta}]=g\#[F^{\eta_{A}}]$ if and
only if $[F^{\Ad g^{-1}\zeta}]=[F^{\eta_{A}}]$.  From
Proposition~\ref{prop:fibreA}, we see that this is the case precisely
when
\begin{align*}
(\lambda\Ad g^{-1}\zeta)(0)&=A\\ \intertext{and} [\Ad
g^{-1}\zeta,A]=0,
\end{align*}
whence the result.
\end{proof}

As an immediate corollary to this result and
Proposition~\ref{prop:ft}, we characterise the maps of finite type in
$\Orbit$:
\begin{cor}\label{cor:fibreC}
$g\#[F^{\eta_{A}}]$ is of finite type if and only if, for some
$d\equiv1\mod k$, there is $\xi\in\Lambda_{d}$ such that
$\xi_{-d}=\Ad g(0)\,A$ and
\[
[\xi,\Ad g\,A]=0.
\]
\end{cor}
For example, the vacuum solutions are themselves of finite type: take
$\xi=\lambda^{-1}A+\lambda\Abar\in\Lambda_{1}$.  These solutions have
{\em constant\/} polynomial Killing fields and so correspond to fixed
points of the Hamiltonian flows on $\Lambda_{1}$.

On the other hand, it is a surprisingly simple matter to construct
solutions in $\Orbit$ which are not of finite type as the following
example shows:
\begin{ex}
Consider the Riemann sphere as a $2$-symmetric space so that
$G=\SU(2)$, $K=S^{1}$ and $B$ consists of diagonal matrices in
$\SL(2,\C)$ of the form
\[
\begin{pmatrix}
r&0\\0&r^{-1}
\end{pmatrix},
\]
where $r>0$.  In this case, $\g_{0}$ consists of trace zero diagonal
matrices and $\g_{1}=\g_{-1}$ of off-diagonal matrices.  We take
\[
A=
\begin{pmatrix}
0&1\\-1&0
\end{pmatrix}
\in\g_{1}
\]
so that $A=\Abar$ and consider the dressing action of
$B\subset\IBloop$ on the vacuum solution $[F^{\eta_{A}}]$.  For $b\in
B$, we know from Corollary~\ref{cor:fibreC} that $b\#[F^{\eta_{A}}]$
is of finite type if and only if we can find $d\equiv1\mod k$ and
$\xi\in\Lambda_{d}$ such that $\xi_{-d}=\Ad b\,A$ and
\[
[\xi,\Ad b\,A]=0.
\]
Writing $\xi=\sum_{|n|\leq d}\lambda^{n}\xi_{n}$, this would force
each $[\xi_{n},\Ad b\,A]$ to vanish and, in particular, since
\[
\overline{\xi_{-d}}=\xi_{d},
\]
we must have
\[
[\overline{\Ad b\, A},\Ad b\, A]=0.
\]
However, when
\[
b=\begin{pmatrix} r&0\\0&r^{-1}
\end{pmatrix}
\]
we have
\[
[\overline{\Ad b\, A},\Ad b\, A]=
\begin{pmatrix}
r^{4}-r^{-4}&0\\0&r^{-4}-r^{4}
\end{pmatrix}
\neq0,
\]
unless $r=1$.  Thus no $b\#[F^{\eta_{A}}]$ is of finite type unless
$b=1$.  This illustrates the complexity of even the simplest case of
the dressing action.
\end{ex}

We are now in a position to identify the stabiliser group $\Stab$:
\begin{thm}\label{thm:stab}
$\Stab$ is the centraliser in $\IBloop$ of
$\lambda^{-1}A\in\initloop$:
\[
\Stab=\{g\in\IBloop\colon\quad \Ad g\,A=A\}.
\]
\end{thm}
\begin{proof}
By Proposition~\ref{prop:fibreB}, $g\in\Stab$ if and only if $\Ad
g(0)\,A=A$ and $[\Ad g\,A,A]$ vanishes.  Thus $\Ad g\, A$ takes
values in $\Ad\GC\,A\cap\ker\ad A$.  On the other hand, since $A$ is
semisimple,
\[
\gc=\ker\ad A\oplus [A,\gc]
\]
from which it follows that $\Ad\GC\,A$ intersects $\ker\ad A$
transversely at $A$ so that $A$ is an isolated point of
$\Ad\GC\,A\cap\ker\ad A$.  The continuity of $\lambda\mapsto\Ad
g(\lambda)\,A$ on the connected set $I_{\epsilon}$ now guarantees
that $\Ad g(\lambda)\,A=A$ for all $\lambda\in I_{\epsilon}$.
\end{proof}

We conclude from this theorem that $\Orbit$ is diffeomorphic to the
adjoint orbit of $\lambda^{-1}A$ in $\initloop$.  In
Section~\ref{sec:flow} we shall make use of this fact to define a
hierarchy of commuting flows on $\Orbit$.

Let us now turn to the main theorem of this section and prove that
all primitive harmonic maps of semisimple finite type lie in some
$\Orbit$.  We begin with a lemma which was proved in the unpublished
thesis of I.~McIntosh.  For completeness of exposition, we shall give
a proof here.
\begin{lem}\label{lem:A-orbit}
Let $X\in\g_{-1}$ be semisimple.  Then there is $A\in\Ad B\,X$ with
$[A,\Abar]=0$.
\end{lem}
\begin{proof}
We begin by finding an element in $\Ad\KC\,X$ with the desired
property.  For this, let $\|\,.\,\|$ denote the norm on $\g_{-1}$
induced by the Killing inner product $(\,.\,)$ on $\g$.  Since $X$ is
semisimple, $\Ad\KC\,X$ is closed in $\g_{-1}$ \cite{Vin76} so that
the restriction of $\|\,.\,\|^{2}$ to $\Ad\KC\, X$ attains a minimum
at some $Y\in\Ad\KC\,X$.  Now, for $\chi\in\sqrt{-1}\k$, we have
\[
0=\left.\frac{d}{dt}\right|_{t=0}\|\Ad\exp t\chi\,Y\|^{2}=
([\chi,Y],\overline{Y})-(Y,[\chi,\overline{Y}])=2(\chi,
[Y,\overline{Y}]).
\]
However, $[Y,\overline{Y}]\in\sqrt{-1}\k$ and the Killing form is
definite there so that $[Y,\overline{Y}]$ vanishes as required.
Thus, for some $k\in\KC$, we have $\Ad k\,X=Y$.  Now put $A=\Ad
k_{K}^{-1}Y$ so that $A=\Ad k_{B}\,X$ and $[A,\Abar]=\Ad
k_{k}^{-1}\,[Y,\overline{Y}]=0$.
\end{proof}

With this in hand, we prove:
\begin{thm}\label{thm:main-orbit}
Let $\eta=\sum_{n\geq
-1}\lambda^{n}\eta_{n}\in\Lambda^{\epsilon'}_{-1,\infty}$ with
$\eta_{-1}$ semisimple.  Then there is $\epsilon\leq\epsilon'$,
$g\in\IBloop$ and a vacuum solution $F^{\eta_{A}}$ such that
\[
g\#[F^{\eta}]=[F^{\eta_{A}}].
\]
\end{thm}
In particular, since
$\lambda^{d-1}\Lambda_{d}\subset\Lambda^{\epsilon'}_{-1,\infty}$ for
all $0<\epsilon'<1$, we have
\begin{cor}\label{cor:main-ft}
Any primitive harmonic map of semisimple finite type is in the
$\IBloop$-orbit of a vacuum solution for some $0<\epsilon<1$.
\end{cor}

Let us turn to the proof of Theorem~\ref{thm:main-orbit}: first, in
view of Lemma~\ref{lem:A-orbit}, we can find $A\in\Ad B\,\eta_{-1}$
such that $[A,\Abar]=0$ and after dressing by an element of $ B$, we
may assume that $\eta_{-1}=A$.  By Proposition~\ref{prop:fibreB}, it
now suffices to find $g\in\IBloop$, some $0<\epsilon\leq\epsilon'$,
such that
\[
\Ad g(0)\,A=A,\quad [A,\Ad g\,\eta]=0.
\]
We shall construct $g$ via the Inverse Function Theorem: since $A$ is
semisimple,
\[
\gc=\ker\ad A\oplus[A,\gc]
\]
and we define $\phi:\ker\ad A\oplus[A,\gc]\to\gc$ by
\[
\phi(x,y)=\Ad\exp(y)\,x.
\]
Observe that $\phi$ is equivariant in the following sense:
\begin{equation}
\label{e:phi-eq}
\omega\tau\phi(x,y)=\phi(\omega\tau x,\tau y),
\end{equation}
for all $(x,y)\in\ker\ad A\oplus[A,\gc]$.

Differentiating $\phi$ at $(A,0)$ gives
\[
\ext_{(A,0)}\phi(v,w)=v+[w,A],
\]
for $(v,w)\in\ker\ad A\oplus[A,\gc]$, so that $\ext_{(A,0)}\phi$ is
an isomorphism.  By the holomorphic Inverse Function Theorem, there
are open neighbourhoods $\Omega_{1}$ of $(A,0)$ and $\Omega_{2}$ of
$A$ such that $\phi:\Omega_{1}\to\Omega_{2}$ is a biholomorphism.
Moreover, since $(A,0)$ is fixed by the order $k$ linear automorphism
$T:(x,y)\to(\omega\tau x,\tau y)$, we may assume, shrinking
$\Omega_{1}$ if necessary, that $\Omega_{1}$ is $T$-stable.

Let $(\psi_{1},\psi_{2})=\phi^{-1}:\Omega_{2}\to\Omega_{1}$ so that,
for $\chi\in\Omega_{2}$,
\[
\chi=\Ad(\exp\psi_{2}(\chi))\psi_{1}(\chi),
\]
or, equivalently,
\begin{equation}
\label{e:psi}
\Ad\exp(-\psi_{2}(\chi))\,\chi=\psi_{1}(\chi)\in\ker\ad A.
\end{equation}
{}From \eqref{e:phi-eq} and the $T$-stability of $\Omega_{1}$ we
deduce that $\psi_{2}$ has the following equivariance property:
\begin{equation}
\label{e:psi-eq}
\psi_{2}(\omega\tau\chi)=\tau\psi_{2}(\chi),
\end{equation}
for all $\chi\in\Omega_{2}$.

Since $\eta\in\Lambda^{\epsilon'}_{-1,\infty}$, $\lambda\eta$ is
holomorphic on $I_{\epsilon'}$ with $(\lambda\eta)(0)=A$ so we can
find $0<\epsilon\leq\epsilon'$ such that $C_{\epsilon}\cup
I_{\epsilon}\subset (\lambda\eta)^{-1}(\Omega_{2})$.  We may
therefore define $g:C_{\epsilon}\cup I_{\epsilon}\to\GC$ by
\[
g(\lambda)=\exp(-\psi_{2}(\lambda\eta(\lambda))).
\]
By construction $g$ is holomorphic on $I_{\epsilon}$ and
$g(0)=\exp(-\psi_{2}(A))=1\in B$ so that
\[
\Ad g(0)\,A=A.
\]
Moreover, from \eqref{e:psi}, for $\lambda\in C_{\epsilon}$, we have
\[
\Ad g(\lambda)\,\eta(\lambda)=
\lambda^{-1}\Ad\exp(-\psi_{2}(\lambda\eta(\lambda)))\,\lambda\eta(\lambda)=
\lambda^{-1}\psi_{1}(\lambda\eta(\lambda))\in\ker\ad A
\]
so that
\[
[A, \Ad g\,\eta]=0
\]
on $C_{\epsilon}$.  Thus $g$ will define our desired element of
$\IBloop$ so long as it satisfies the equivariance condition
$g(\omega\lambda)=\tau g(\lambda)$.  For this, recall that
$\eta(\omega\lambda)=\tau\eta(\lambda)$ so that, using
\eqref{e:psi-eq},
\begin{align*}
g(\omega\lambda)&=\exp(-\psi_{2}(\omega\lambda\eta(\omega\lambda)))=
\exp(-\psi_{2}(\omega\tau\lambda\eta(\lambda)))\\
&=\exp(-\tau\psi_{2}(\lambda\eta(\lambda)))=\tau g(\lambda)
\end{align*}
as required.  This completes the proof of
Theorem~\ref{thm:main-orbit}.

\begin{rem}
Corollary~\ref{cor:main-ft} is an extension of results of
Dorfmeister--Wu \cite{DorWu93} and Wu \cite{Wu93}.  In those papers,
a similar result was proved by different methods for $G=\SU(2)$ and
$G=\SO(5)$, respectively, with $k$-symmetric structure given by the
Coxeter--Killing automorphism.  Geometrically, these cases correspond
to non-conformal harmonic maps into $S^{2}$ (i.e. Gauss maps of
constant mean curvature surfaces in $\R^{3}$) and (twistor lifts of)
minimal non-superminimal maps into $S^{4}$.  In view of
\cite{Bur93,BurFerPed93}, the present result accounts for all
harmonic non-isotropic $2$-tori in any sphere or complex projective
space.
\end{rem}

To summarise: our results give a rather complete picture of the
dressing orbits through primitive harmonic maps $[F^{\eta}]$ with
$(\lambda\eta)(0)$ semisimple.  Any non-constant map of this type
lies in some $\Orbit$ and each such orbit is infinite-dimensional.
Moreover, if $A$ and $A'$ both generate vacuum solutions and lie in
the same $\Ad B$-orbit, we clearly have
\[
\Orbit=\cal{O}_{A'}.
\]
Finally, for $w\in\C^{*}$, we can scale the parameter $z$ on $\R^{2}$
to identify $\Orbit$ and $\cal{O}_{wA}$.  We therefore conclude that
the set of essentially different dressing orbits of vacuum solutions
is parametrised by the projective $\Ad B$-orbits of semisimple
elements of $\g_{1}$.  This is a finite-dimensional family.

We conclude this section with an application of
Theorem~\ref{thm:main-orbit} which provides a simplified version of
Lemma~\ref{lem:fibre}.  Recall that a semisimple $A\in\gc$ is {\em
regular\/} if its centraliser $\ker\ad A$ is abelian.
\begin{prop}
Let $\zeta$, $\eta\in\initloop$ with
$(\lambda\zeta)(0)=(\lambda\eta)(0)$ regular semisimple.  Then
$[F^{\zeta}]=[F^{\eta}]$ if and only if $[\zeta,\eta]=0$.
\end{prop}
\begin{proof}
In view of Lemma~\ref{lem:fibre}, only the forward implication
requires proof.  So suppose that $[F^{\zeta}]=[F^{\eta}]$.  By
Theorem~\ref{thm:main-orbit}, there is a vacuum solution
$[F^{\eta_{A}}]$ and $g\in\Lambda^{\delta}_{I}G^{C}_{\tau}$, for
$0<\delta\leq\epsilon$, such that
\[
[F^{\zeta}]=[F^{\eta}]=g\#[F^{\eta_{A}}].
\]
Thus, by Proposition~\ref{prop:fibreB}, $(\lambda\zeta)(0)=\Ad
g(0)\,A$, whence $A$ is regular semisimple, and
\[
[\zeta,\Ad g\,A]=[\eta,\Ad g\,A]=0.
\]
Thus $\Ad g^{-1}\,\zeta$ and $\Ad g^{-1}\,\eta$ take values in
$\ker\ad A$ and so commute.  Thus $[\zeta,\eta]=0$ as required.
\end{proof}
\begin{ex}
When $\tau$ is the Coxeter--Killing automorphism of a simple group
$G$, then any semisimple element in $\g_{-1}$ is regular
\cite{Kos59}.
\end{ex}

\section{Higher flows}
\label{sec:flow}

It is characteristic of ``completely integrable'' systems of partial
differential equations that the solution set admits a hierarchy of
commuting flows.  In this section, we shall define such a hierarchy
on each dressing orbit $\Orbit$ and characterise the primitive
harmonic maps of finite type in $\Orbit$ as precisely those whose
orbit under these flows is finite-dimensional.

\begin{notation}
We denote by $[g]$ the coset of $g$ in $\IBloop/\Stab$.
\end{notation}

Recall from Section~\ref{sec:vac} that we have a diffeomorphism
$\IBloop/\Stab\cong\Orbit$ given by
\[
[g]\mapsto g\#[F^{\eta_{A}}].
\]
In view of the identification $\IBloop\cong\Eloop\backslash\Loop$, we
have a right action of $\Loop$ on $\IBloop$ given by
\[
(g,h)\mapsto (gh)_{I}.
\]
We now identify an abelian subgroup of $\Loop$ for which this action
descends to one on $\Orbit$.  For this, let $\centr\subset\gc$ denote
the centre of the centraliser of $A$ and define $\comm\subset\gloop$
by
\[
\comm=\{\zeta\in\gloop\colon\text{for some $N\in\N$, $\zeta(\lambda)=
\sum_{n\geq-N}\lambda^{n}\zeta_{n}$, for $\lambda\in C_{\epsilon}$,
with $\zeta_{n}\in\centr$ for all $n$.}\}
\]
Moreover, let $\Comm\subset\Loop$ be the abelian subgroup obtained by
exponentiating $\comm$.

Since $A$ is semisimple, the centraliser of $A$ in $\Ad\GC$ is
connected \cite[Lemma~5]{Kos63} from which it follows that $\Comm$
commutes with $\Stab$.  Thus, for $\exp\zeta\in\Comm$ and
$[g]=[g']\in\Orbit$, we have $g=g'\gamma$, some $\gamma\in\Stab$, so
that
\[
(g\exp\zeta)_{I}=(g'\gamma\exp\zeta)_{I}=(g'\exp(\zeta)\gamma)_{I}=
(g'\exp\zeta)_{I}\gamma
\]
whence
\[
[(g\exp\zeta)_{I}]=[(g'\exp\zeta)_{I}]
\]
and we have proved:
\begin{prop}
$\Comm$ acts on $\Orbit\cong\IBloop/\Stab$ by
\[
\exp\zeta\cdot[g]=[(g\exp\zeta)_{I}].
\]
\end{prop}

\begin{rem}
It is easy to see that $Z\cap\Stab=Z\cap\IBloop$ acts trivially on
$\Orbit$.  With a little more work, one can also show that any $z\in
Z\setminus\Iloop$ acts non-trivially on $\Orbit$.  In particular,
when $\centr\cap\g_{0}=\{0\}$ (which is the case, for instance, for
abelian Toda fields), we deduce that there is an effective action of
$Z/Z\cap\IBloop$ on $\Orbit$.
\end{rem}
Our main result concerning the action of $\Comm$ is the following
characterisation of the primitive harmonic maps of finite type in
$\Orbit$:
\begin{thm}\label{thm:comm-main}
$[F]\in\Orbit$ is of finite type if and only if the $\Comm$-orbit of
$[F]$ is finite-dimensional.
\end{thm}
\begin{rem}
In their study of the sinh-Gordon equation, Dorfmeister--Wu prove
this result by different methods for the case $G=\SU(2)$ equipped
with the Coxeter--Killing automorphism.
\end{rem}

We break up our proof of Theorem~\ref{thm:comm-main} into a sequence
of steps.  We begin with a simple technical lemma:
\begin{lem}\label{lem:E}
Let $d\geq 1\in\N$ and suppose that $\zeta\in\gloop$ with
$\lambda^{d}\zeta$ holomorphic on $I_{\epsilon}$.  Then
$\zeta_{E}\in\Lambda_{d}$ and
\[
(\zeta_{E})_{-d}=(\lambda^{d}\zeta)(0).
\]
\end{lem}
\begin{proof}
On $C_{\epsilon}$, we can write
$\zeta=\sum_{n\geq-d}\lambda^{n}\zeta_{n}$ and it is easy to see that
\[
\zeta_{E}=\sum_{-d\leq n\leq-1}\lambda^{n}\zeta_{n}+(\zeta_{0})_{\k}+
\sum_{1\leq n\leq d}\lambda^{n}\overline{\zeta_{-n}}
\]
from which the lemma follows immediately.
\end{proof}

We can now prove one part of the theorem:
\begin{prop}\label{prop:comm-easy}
If $[F]=g\#[F^{\eta_{A}}]\in\Orbit$ has finite-dimensional
$\Comm$-orbit then $[F]$ is of finite type.
\end{prop}
\begin{proof}
According to Corollary~\ref{cor:fibreC}, it suffices to find
$d\equiv1\mod k$ and $\xi\in\Lambda_{d}$ such that $\xi_{-d}=\Ad
g(0)\,A$ and
\[
[\xi,\Ad g\,A]=0.
\]
By hypothesis, the map $\comm\to T_{[F]}\Orbit$ given by
\[
\zeta\mapsto\left.\frac{d}{dt}\right|_{t=0}\exp t\zeta\cdot [g]
\]
has finite-dimensional image.  We use left translation by $g$ to
identify $T_{[F]}\Orbit$ with
$T_{[F^{\eta_{A}}]}\Orbit\cong\ibloop/\stab$ (here $\stab$ is the Lie
algebra of $\Stab$) and then an easy calculation shows that the
resulting map $\comm\to\ibloop/\stab$ is given by
\begin{equation}
\label{e:comm-stab}
\zeta\mapsto\Ad g^{-1}\,(\Ad g\,\zeta)_{I}\mod\stab.
\end{equation}
Let $\cal{N}\subset\comm$ denote the kernel of the linear map
\eqref{e:comm-stab}.  Since $\cal N$ has finite codimension, there is
a monic polynomial $p$ such that
$\zeta=p(\lambda^{-k})\lambda^{-1}A\in\comm$ lies in $\cal N$ so that
\[
\Ad g^{-1}\,(\Ad g\,\zeta)_{I}\in\stab
\]
or, equivalently,
\[
[\Ad g^{-1}\,(\Ad g\,\zeta)_{I},A]=0.
\]
We therefore see that $[(\Ad g\,\zeta)_{I},\Ad g\,A]=0$.  On the
other hand, since $\zeta\in\comm$, we have $[\zeta,A]=0$ and we can
conclude that
\begin{equation}
\label{e:E-A}
[(\Ad g\,\zeta)_{E},\Ad g\,A]=0.
\end{equation}
Now set $d=k\deg p+1\equiv1\mod k$ and observe that
$\lambda^{d}\zeta$ is holomorphic (indeed polynomial) on
$I_{\epsilon}$ with $(\lambda^{d}\zeta)(0)=A$.  It follows that
$\lambda^{d}\Ad g\,\zeta$ is holomorphic on $I_{\epsilon}$ with
$(\lambda^{d}\Ad g\,\zeta)(0)=\Ad g(0)\,A$ so that, from
Lemma~\ref{lem:E}, we see that $\xi=(\Ad g\,\zeta)_{E}\in\Lambda_{d}$
with $\xi_{-d}=\Ad g(0)\,A$.  This, taken together with
\eqref{e:E-A}, establishes the proposition.
\end{proof}
\begin{rem}
In case that $G=\SO(5)$ equipped with the Coxeter--Killing
automorphism, this proposition was proved by Wu \cite{Wu93} using
different methods.
\end{rem}

It remains to prove the converse of Proposition~\ref{prop:comm-easy}.
For this, we use an argument inspired by ideas of McIntosh
\cite{McI94B}.  Fix $d\equiv1\mod k$ and let $\dOrbit$ consist of
those maps in $\Orbit$ which admit a polynomial Killing field
$\R^{2}\to\Lambda_{d}$.  $\dOrbit$ is the image under
$\initloop\to\Harm$ of a subvariety of $\lambda^{d-1}\Lambda_{d}$ and
so is finite-dimensional.  With this in mind, the proof of
Theorem~\ref{thm:comm-main} is completed by the following

\begin{prop}
$\dOrbit$ is $\Comm$-stable.
\end{prop}
\begin{proof}
Let $[g]\in\dOrbit$ so that there is $\xi\in\Lambda_{d}$ such that
$(\lambda^{d}\xi)(0)=\Ad g(0)\,A$ and
\[
[\xi,\Ad g\,A]=0.
\]
Let $\zeta\in\comm$ and set
\[
\hat{\xi}=\Ad(g\exp\zeta)_{I}\Ad g^{-1}\,\xi= \Ad(\exp\Ad
g\,\zeta)_{I}\,\xi.
\]
Clearly we have
\[
[\hat{\xi},\Ad(g\exp\zeta)_{I}\,A]=[\Ad g^{-1}\,\xi,A]=0.
\]
It therefore suffices to show that $\hat{\xi}\in\Lambda_{d}$ and that
$\hat{\xi}_{-d}=\Ad(g\exp\zeta)_{I}(0)\,A$ to conclude that
$\exp\zeta\cdot[g]\in\dOrbit$.

For this, first note that since $\Ad g^{-1}\,\xi$ centralises $A$,

$[\zeta,\Ad g^{-1}\,\xi]=0$ whence
\[
\Ad(g\exp\zeta)\,\Ad g^{-1}\,\xi=\Ad g\,\Ad g^{-1}\,\xi=\xi
\]
so that
\[
\Ad(\exp\Ad g\,\zeta)_{E}\,\Ad(\exp\Ad g\,\zeta)_{I}\,\xi=\xi
\]
from which we conclude that
\[
\hat{\xi}=\Ad(\exp\Ad g\,\zeta)_{I}\,\xi= \Ad(\exp\Ad
g\,\zeta)^{-1}_{E}\,\xi\in\eloop.
\]
On the other hand,
\[
\lambda^{d}\hat{\xi}=\Ad(g\exp\zeta)_{I}\,\Ad g^{-1}\,\lambda^{d}\xi
\]
which is holomorphic on $I_{\epsilon}$ and
\[
(\lambda^{d}\hat{\xi})(0)=\Ad(g\exp\zeta)_{I}(0)\,\Ad
g^{-1}(0)\,\xi_{-d}= \Ad(g\exp\zeta)_{I}(0)\,A.
\]
It now follows from Lemma~\ref{lem:E} that $\hat{\xi}$ has the
required properties and the proof is complete.
\end{proof}

\section*{Appendix: nilpotency of $\boldsymbol{\eta_{-1}}$ and finite uniton
number} \setcounter{equation}{0} \setcounter{thm}{0}
\renewcommand{\thesection}{A}

In Section~\ref{sec:vac}, we saw that the primitive harmonic maps
$[F^{\eta}]$ with $(\lambda\eta)(0)$ semisimple comprise a
finite-dimensional family of infinite-dimensional dressing orbits.
By contrast, we now consider the case where $(\lambda\eta)(0)$ is
nilpotent (which is necessarily the case if the harmonic map
$[F^{\eta}]$ has finite energy).  Here our results are less complete
and somewhat confused but examples and partial results indicate that
the picture is completely different.

\begin{ex}
Let $\eta=\lambda^{-1}\eta_{-1}\in\initloop$ with $\eta_{-1}$
nilpotent so that $(\ad\eta_{-1})^{\ell}=0$, say.  Comparing
coefficients in \eqref{e:comm}, we conclude from
Lemma~\ref{lem:fibre} that $[F^{\zeta}]=[F^{\eta}]$ if and only if
$\zeta_{-1}=\eta_{-1}$ and
\[
(\ad\eta_{-1})^{n}\zeta_{n-1}=0,
\]
for $n<\ell$.  From this it is easy to see that
$g\#[F^{\eta}]=[F^{\eta}]$ for all $g\in\IBloop$ whose $\ell$-jet at
zero coincides with that of $1$.  Such $g$ comprise a normal subgroup
of $\IBloop$ of finite codimension so that the action of $\IBloop$ on
$[F^{\eta}]$ reduces to that of a finite-dimensional quotient group.
\end{ex}

This is reminiscent of the behaviour of harmonic maps with finite
uniton number in the sense of Uhlenbeck \cite{Uhl89}.  A primitive
harmonic map has finite uniton number if it has an extended framing
$F$ which can be written as a Laurent polynomial of fixed degree in
$\lambda$:
\[
F(z)=\sum_{|n|\leq d}\lambda^{n}F_{n}(z),
\]
where we fix a faithful representation of $\GC$ to make sense of the
Laurent expansion.  We call the minimal such $d$ the {\em uniton
number\/} of $[F]$.

Concerning these maps, Uhlenbeck \cite{Uhl89} and Bergveldt--Guest
\cite{BerGue91} prove
\begin{prop}
The action of $\Iloop$ preserves uniton number and the action on the
space of maps of fixed uniton number reduces to that of a finite
dimensional quotient group.
\end{prop}

Since the orbits $\Orbit$ are infinite-dimensional, we deduce the
following result announced in \cite{BurPed94}:
\begin{thm}
No $[F^{\eta}]$ with $(\lambda\eta)(0)$ non-zero semisimple and, in
particular, no non-constant primitive harmonic map of semisimple
finite type, has finite uniton number.
\end{thm}

By contrast, we have the following model result which we state for
$G=\SU(2)$ for simplicity of exposition.
\begin{prop}\label{prop:f-u}
Let $[F^{\eta}]:\R^{2}\to S^{2}$ be harmonic with $\eta_{-1}$
nilpotent.  Then $[F^{\eta}]$ has finite uniton number.
\end{prop}
\begin{proof}
Set $F=F^{\eta}$.  In this setting the reality condition reads
\[
F_{1/\bar{\lambda}}^{*}=F^{-1}_{\lambda}
\]
from which it is easy to see that $F$ is a Laurent polynomial
precisely when, for some $d$, both $\lambda^{d}F_{\lambda}$ and
$\lambda^{d}F^{-1}_{\lambda}$ are holomorphic (as matrix valued
functions) at $\lambda=0$.  In view of the definition of $F$, this
amounts to demanding that $\lambda^{d}\exp(\pm z\eta)$ be holomorphic
at zero.

To prove this, we first note that after dressing with a constant
element of $\KC$, we may assume that
\[
\eta_{-1}=
\begin{pmatrix}
0&1\\0&0
\end{pmatrix}
\quad\text{or}\quad
\begin{pmatrix}
0&0\\1&0
\end{pmatrix}.
\]
For definiteness, we assume
\[
\eta_{-1}=
\begin{pmatrix}
0&1\\0&0
\end{pmatrix}
\]
and write
\[
\eta(\lambda)=
\begin{pmatrix}
a(\lambda)&\lambda^{-1}+b(\lambda)\\ c(\lambda)&-a(\lambda)
\end{pmatrix}
\]
where $a$, $b$, $c$ are holomorphic on $I_{\epsilon}$.  The
equivariance condition on $\eta$ means that $a$ is an even function
while $b$ and $c$ are odd.  In particular, $c(0)=0$.

Now set
\[
\gamma(\lambda)=
\begin{pmatrix}
1&0\\0&\lambda^{-1}
\end{pmatrix}
\]
and observe that
\[
\Ad\gamma(\lambda)\,\eta(\lambda)=
\begin{pmatrix}
a(\lambda)&1+\lambda b(\lambda)\\ \lambda^{-1}c(\lambda)&-a(\lambda)
\end{pmatrix}
\]
is holomorphic at $\lambda=0$ since $\lambda^{-1}c(\lambda)$ is.
Thus
\[
\gamma(\lambda)\exp(\pm z\eta(\lambda))\gamma(\lambda)^{-1}
\]
is holomorphic at zero whence $\lambda\exp(\pm z\eta(\lambda))$ is
also.
\end{proof}
\begin{rem}
The trick of conjugating with $\gamma$ amounts to passing between the
principal and standard realisations of $\Lambda\SL(2,\C)$ in the
sense of Wilson \cite{Wil85}.
\end{rem}

It is not difficult to extend the above argument to numerous other
cases which include that where $G/K$ is a full flag manifold with
Coxeter--Killing automorphism and $\eta_{-1}$ is principal nilpotent.
In view of this, one might feel tempted (as we were!) to conjecture
that any primitive harmonic map $[F^{\eta}]$ with $\eta_{-1}$
nilpotent has finite uniton number.  However, this is not the case as
the following analysis of harmonic maps $\R^{2}\to\SU(2)$ shows.

Set $K=\SU(2)$ and put $G=K\times K$.  We view $K$ as the Riemannian
symmetric space $G/K$ with involution $\tau:G\to G$ given by
$(k_{1},k_{2})\mapsto (k_{2},k_{1})$ and coset projection $G\to K$
given by $(k_{1},k_{2})\mapsto k_{1}k_{2}^{-1}$.

Define $\Lambda_{\mathrm{hol}}K$ by
\[
\Lambda_{\mathrm{hol}}K=\{k:\C^{*}\to\KC\colon\text{$k$ is
holomorphic and $\overline{k(\lambda)}=k(1/\bar{\lambda})$}\}.
\]
It is easy to see that $\Hloop$ is given by
\[
\Hloop=\{(k(-\lambda),k(\lambda))\colon\quad
k\in\Lambda_{\mathrm{hol}}K\}
\]
so that projection onto the second factor gives an isomorphism
$\Hloop\cong\Lambda_{\mathrm{hol}}K$.  Similarly, we can define
subgroups $\Lambda^{\epsilon}_{E}K$ and so on, just by dropping the
equivariance conditions on the loops $C^{(\epsilon)}\to\KC$ and thus
obtain isomorphisms $\Eloop\cong\Lambda^{\epsilon}_{E}K$ and so on.
Under these isomorphisms, the harmonic map $\phi:M\to K$
corresponding to the extended framing $F:M\to\Lambda_{\mathrm{hol}}K$
is given by
\[
\phi=F_{-1}F_{1}^{-1}
\]
(and, indeed, the {\em extended solution\/} corresponding to $\phi$
in the sense of Uhlenbeck is the map into the based loop group
$\Omega K$ given by $F_{\lambda}F_{1}^{-1}:M\to\Omega K$).

Having adopted this view-point, we may think of $\eta\in\initloop$ as
a holomorphic map $\eta:I_{\epsilon}\setminus\{0\}\to\kc$ with at
most a simple pole at zero.  We now have
\begin{prop}\label{prop:holo-det}
$[F^{\eta}]$ has finite uniton number if and only if $\det\eta$ is
holomorphic on $I_{\epsilon}$.
\end{prop}
\begin{proof}
Arguing as in Proposition~\ref{prop:f-u}, we see that $[F^{\eta}]$
has finite uniton number if and only if, for some $d\in\N$,
$\lambda^{d}\exp(z\eta)$ is holomorphic on $I_{\epsilon}$, for all
$z\in\C$.

For $A\in\Sl(2,\C)$, one has
\[
A^{2}=(-\det A)\id
\]
from which we deduce, setting $D=\sqrt{-\det A}$, that
\[
\exp A=\cosh(D)\id + \frac{\sinh(D)}{D}A.
\]
Thus
\begin{equation}
\label{e:expA}
\exp(z\eta)=\cosh(\sqrt{-\det z\eta})\id+\frac{\sinh(\sqrt{-\det
z\eta})}{\sqrt{-\det z\eta}}z\eta.
\end{equation}
If $\det\eta$ is holomorphic on $I_{\epsilon}$, we see that
$\lambda\exp(z\eta)$ is also since, by hypothesis, $\lambda\eta$ is
holomorphic on $I_{\epsilon}$.  Thus, in this case, $[F^{\eta}]$ has
finite uniton number.

Conversely, if $\det\eta$ has a pole at zero, it is easy to see from
\eqref{e:expA} that at least one matrix entry of $\exp(z\eta)$ has an
essential singularity at zero so that $[F^{\eta}]$ has infinite
uniton number.
\end{proof}

With this in hand, it is easy to find $[F^{\eta}]$ with $\eta_{-1}$
nilpotent and infinite uniton number:
\begin{ex}
Define $\eta\in\initloop$ by
\[
\eta(\lambda)=\begin{pmatrix} 0&\lambda^{-1}\\ 1&0
\end{pmatrix}
\]
so that
\[
\eta_{-1}=\begin{pmatrix} 0&1\\0&0
\end{pmatrix}
\]
is nilpotent while $\det\eta=-\lambda^{-1}$ whence $[F^{\eta}]$ has
infinite uniton number by Proposition~\ref{prop:f-u}.
\end{ex}
\ifx\undefined\bysame
\newcommand{\bysame}{\leavevmode\hbox to3em{\hrulefill}\,}
\fi


\begin{thebibliography}{10}

\bibitem{BerGue91}
M.J. Bergveldt and M.A. Guest, {\em {A}ction of loop groups on harmonic maps},
  Trans.\ Amer.\ Math.\ Soc. {\bf 326} (1991), 861--886.

\bibitem{BolPedWoo93}
J.~Bolton, F.~Pedit, and L.M. Woodward, {\em {M}inimal surfaces and the affine
  {T}oda field model}, J. reine angew.\ Math. (to appear).

\bibitem{Bur93}
F.E. Burstall, {\em {H}armonic tori in spheres and complex projective spaces},
  preprint, 1993.

\bibitem{BurFerPed93}
F.E. Burstall, D.~Ferus, F.~Pedit, and U.~Pinkall, {\em {H}armonic tori in
  symmetric spaces and commuting {H}amiltonian systems on loop algebras}, Ann.\
  of Math. {\bf 138} (1993), 173--212.

\bibitem{BurPed94}
F.E. Burstall and F.~Pedit, {\em {H}armonic maps via
  {A}dler--{K}ostant--{S}ymes theory}, {H}armonic maps and integrable systems
  (A.P. Fordy and J.C. Wood, eds.), Aspects of Math., vol.~23, Vieweg,
  Braunschweig, Wiesbaden, 1994, pp.~221--272.

\bibitem{DorPedWu94}
J.~Dorfmeister, F.~Pedit, and H.~Wu, {\em {W}eierstrass type representation of
  harmonic maps into symmetric spaces}, GANG preprint III.25, 1994.

\bibitem{DorWu93}
J.~Dorfmeister and H.~Wu, {\em {C}onstant mean curvature surfaces and loop
  groups}, J. reine angew.\ Math. {\bf 440} (1993), 43--76.

\bibitem{FerPedPin92}
D.~Ferus, F.~Pedit, U.~Pinkall, and I.~Sterling, {\em {M}inimal tori in
  {${S}^4$}}, J. reine angew.\ Math. {\bf 429} (1992), 1--47.

\bibitem{GueOhn93}
M.A. Guest and Y.~Ohnita, {\em {G}roup actions and deformations for harmonic
  maps}, J. Math.\ Soc.\ Japan {\bf 45} (1993), 671--704.

\bibitem{Hit90}
N.J. Hitchin, {\em {H}armonic maps from a {$2$}-torus to the {$3$}-sphere}, J.
  Diff.\ Geom. {\bf 31} (1990), 627--710.

\bibitem{JenLia94}
G.R. Jensen and R.~Liao, {\em {F}amilies of flat minimal tori in {${\bf
  {C}}{P}^n$}}, Preprint, 1994.

\bibitem{Kos59}
B.~Kostant, {\em {T}he principal three dimensional subgroup and the {B}etti
  numbers of a complex simple {L}ie group}, Amer.\ J. Math. {\bf 81} (1959),
  973--1032.

\bibitem{Kos63}
\bysame, {\em {L}ie group representations on polynomial rings}, Amer.\ J. Math.
  {\bf 86} (1963), 327--402.

\bibitem{Kow80}
O.~Kowalski, {\em {G}eneralized symmetric spaces}, Lect.\ Notes in Math., vol.
  805, Springer, Berlin, Heidelberg, New York, 1980.

\bibitem{LuWei90}
J.-H. Lu and A.~Weinstein, {\em {P}oisson {L}ie groups, dressing
  transformations and {B}ruhat decompositions}, J. Diff.\ Geom. {\bf 31}
  (1990), 501--526.

\bibitem{McI94}
I.~McIntosh, {\em {G}lobal solutions of the elliptic 2{D} periodic {T}oda
  lattice}, Nonlinearity {\bf 7} (1994), 85--108.

\bibitem{McI94B}
\bysame, {\em {O}n the existence of superconformal {$2$}-tori in complex
  projective space}, Bath preprint, 1994.

\bibitem{PinSte89}
U.~Pinkall and I.~Sterling, {\em {O}n the classification of constant mean
  curvature tori}, Ann.\ of Math. {\bf 130} (1989), 407--451.

\bibitem{Poh76}
K.~Pohlmeyer, {\em {I}ntegrable {H}amiltonian systems and interactions through
  quadratic constraints}, Commun.\ Math.\ Phys. {\bf 46} (1976), 207--221.

\bibitem{PreSeg86}
A.N. Pressley and G.~Segal, {\em {L}oop {G}roups}, Oxford Math.\ Monographs,
  Clarendon Press, Oxford, 1986.

\bibitem{Sym80}
W.~Symes, {\em {S}ystems of {T}oda type, inverse spectral problems and
  representation theory}, Invent.\ Math. {\bf 159} (1980), 13--51.

\bibitem{Uda94}
S.~Udagawa, {\em {H}armonic maps from a two-torus into a complex {G}rassmann
  manifold}, preprint, 1994.

\bibitem{Uda94a}
\bysame, {\em {H}armonic tori in quaternionic projective {$3$}-spaces},
  preprint, 1994.

\bibitem{Uhl89}
K.~Uhlenbeck, {\em {H}armonic maps into {L}ie groups (classical solutions of
  the chiral model)}, J. Diff.\ Geom. {\bf 30} (1989), 1--50.

\bibitem{Vin76}
E.B. Vinberg, {\em {T}he {W}eyl group of a graded {L}ie algebra}, Izv.\ Akad.\
  Nauk SSSR Ser.\ Mat. {\bf 40} (1976), 488--526, 709, Math.\ USSR
  Izv. {\bf 10} (1976), 463--495.

\bibitem{Wil85}
G.~Wilson, {\em {I}nfinite-dimensional {L}ie groups and algebraic geometry in
  soliton theory}, Phil.\ Trans.\ R. Soc.\ Lond.\ A {\bf 315} (1985), 393--404.

\bibitem{Wu93}
H.~Wu, {\em {B}anach manifolds of minimal surfaces in the {$4$}-sphere},
  {D}ifferential {G}eometry (R.E. Greene and S.T. Yau, eds.), Proc.\ Symp.\
  Pure Math., vol.~54, Amer.\ Math.\ Soc., Providence RI, 1993, pp.~513--539.

\bibitem{ZakSha78}
V.E. Zakharov and A.B. Shabat, {\em {I}ntegration of non-linear equations of
  mathematical physics by the method of inverse scattering, {I}{I}}, Funkts.\
  Anal.\ i Prilozhen {\bf 13} (1978), 13--22, Funct.\ Anal.\ Appl.
  {\bf 13} (1979), 166--174.

\end{thebibliography}
\end{document}